\documentclass[10pt,twocolumn,twoside]{IEEEtran}
\IEEEoverridecommandlockouts                            
\usepackage[english]{babel}
\usepackage[dvipsnames]{xcolor}
\usepackage{thmtools}
\usepackage{graphicx}
\usepackage[normalem]{ulem}
\usepackage{amsmath}
\usepackage{amssymb}
\usepackage{hyperref}
\usepackage{amsfonts}
\usepackage{bbm}
\usepackage{fge}
\usepackage{enumerate}
\usepackage{euscript}
\usepackage{empheq}
\usepackage{mathtools}
\usepackage{tikz}
\usepackage{pgfplots}
\usepackage{hyperref}
\usepackage{balance}
\usepackage{cite}
\usepackage{booktabs}
\usepackage{xfrac}
\usepackage{algorithm}
\usepackage{multirow}
\usepackage[noend]{algpseudocode}
\usetikzlibrary{arrows}
\usetikzlibrary{shapes}
\usetikzlibrary{patterns}
\usepgfplotslibrary{fillbetween}
\usetikzlibrary{arrows.meta}
\usepgfplotslibrary{colorbrewer}
\usepackage{framed}
\usepackage{dblfloatfix}
\usepackage[makeroom]{cancel}
\usepackage{xcolor}
\usepackage{tikz}
\usepackage{graphicx}
\usepackage{pgfplots}
\pgfplotsset{compat=1.16}
\usepackage[bb=boondox]{mathalfa}
\usepackage{nicefrac}
\usepackage{upgreek}
\usepackage{relsize}
\usepackage{cuted}
\usepackage{pgfplotstable}

\usepackage{nomencl}
\makenomenclature
\usepackage{etoolbox}
\renewcommand\nomgroup[1]{%
  \item[\itshape
  \ifstrequal{#1}{A}{Sets}{%
  \ifstrequal{#1}{B}{Parameters}{%
  \ifstrequal{#1}{C}{Variables}{}}}%
]}
\makeatletter
\patchcmd{\thenomenclature}{\section*}{\section}{}{}
\makeatother

\newtheorem{theorem}{Theorem}

\newtheorem{remark}{Remark}
\newtheorem{lemma}{Lemma}
\newtheorem{corollary}{Corollary}

\allowdisplaybreaks

\newcommand{\minimize}[1]{\underset{{#1}}{\text{min}}}
\newcommand{\maximize}[1]{\underset{{#1}}{\text{max}}}

\newcommand{\col}[1]{\langle #1 \rangle}
\newcommand{\mysum}[2]{\textstyle\sum_{#1}^{#2}}

\newcommand{\braceit}[1]{\left({#1}\right)}
\newcommand{\st}{\text{s.t.}}

\usepackage{mathtools}
\newcommand\norm[1]{\lVert#1\rVert}

\usepackage{esvect}

\makeatletter
\newcommand{\oset}[3][0ex]{%
  \mathrel{\mathop{#3}\limits^{
    \vbox to#1{\kern-2\ex@
    \hbox{$\scriptstyle#2$}\vss}}}}
\makeatother

\usetikzlibrary{decorations.pathreplacing}
\makeatletter
\def\underbrace#1{%
   \@ifnextchar_{\tikz@@underbrace{#1}}{\tikz@@underbrace{#1}_{}}}
\def\tikz@@underbrace#1_#2{%
   \tikz[baseline=(a.base)] {\node[inner sep=2] (a) {\(#1\)};
   \draw[line cap=round,decorate,decoration={brace,amplitude=3pt}]
     (a.south east) -- node[pos=0.5,below,inner sep=5pt] {\(\scriptstyle #2\)} (a.south west);}}

\definecolor{backg_vp_plot}{RGB}{237,243,254}
\definecolor{red}{HTML}{e05a87} 
\definecolor{maincolor}{HTML}{032F99}

\usepgfplotslibrary{colorbrewer}
\pgfplotsset{compat = 1.15, cycle list/Set1-8} 
\usetikzlibrary{pgfplots.statistics} 
\usepackage{pgfplotstable}
\usepackage{filecontents}

\title{\LARGE \bf 
Stochastic Control and Pricing for Natural Gas Networks
}

\author{Vladimir Dvorkin, Anubhav Ratha, Pierre Pinson and Jalal Kazempour% <-this % stops a space
\vspace{-0.5cm}
 \thanks{%Manuscript received September 30, 2020. 
 The authors would like to thank Ana Virag and H\'{e}l\`{e}ne Le Cadre for their helpful comments, Jochen Stiasny and Tue Jensen for the advice on visualizations, Roberth Mieth and Yury Dvorkin for discussions.}
\thanks{V. Dvorkin is with MIT Energy Initiative and Laboratory for Information and Decision Systems, Cambridge, MA, USA. Email: {\tt dvorkin@mit.edu}}
\thanks{A. Ratha and J. Kazempour are with the Department of Electrical Engineering, Technical University of Denmark, Lyngby, Denmark.
         Email: {\tt \{arath,seykaz\}@elektro.dtu.dk}}%
\thanks{P. Pinson is with the Department of Technology, Management and Economics, Technical University of Denmark, Lyngby, Denmark. Email: {\tt ppin@dtu.dk}}
\thanks{A. Ratha is also with the Flemish Institute of Technological Research (VITO), Boeretang 200, 2400 Mol, Belgium and with EnergyVille, ThorPark 8310, 3600 Genk, Belgium.
        {\tt \{anubhav.ratha\}@vito.be}}%
}

\begin{document}
\begingroup
\allowdisplaybreaks

\maketitle

\begin{abstract}
We propose stochastic control policies to cope with uncertain and variable gas extractions in natural gas networks. Given historical gas extraction data, these policies are optimized to produce the real-time control inputs for nodal gas injections and for pressure regulation rates by compressors and valves. We describe the random network state as a function of control inputs, which enables a chance-constrained optimization of these policies for arbitrary network topologies. This optimization ensures the real-time gas flow feasibility and a minimal variation in the network state up to specified feasibility and variance criteria. Furthermore, the chance-constrained optimization provides the foundation of a stochastic pricing scheme for natural gas networks, which improves on a deterministic market settlement by offering the compensations to network assets for their contribution to uncertainty and variance control. We analyze the economic properties, including efficiency, revenue adequacy and cost recovery, of the proposed pricing scheme and make them conditioned on the network design. 
\end{abstract}

\begin{IEEEkeywords}
Chance-constrained programming, conic duality, gas pricing, natural gas network, uncertainty, variance.
\end{IEEEkeywords}

%%%%%%%%%%%%%%%%%%%%%%%%%%%%%%%%%%%%%%%%%%%%%%%%%%%%%%%%%%%%%%%%%%%%%%%%%%%%%%%%
\section{Introduction} \label{sec:Introduction}

Deterministic operational and market-clearing practices of the natural gas network operators struggle with the growing uncertainty and variability of natural gas extractions \cite{outlook20192019}. Ignorance of the uncertain and variable extractions results in technical and economical failures, as demonstrated by the congested network during the 2014 polar vortex event in the United States \cite{interconnection2014analysis}. The recent study \cite{bent2018joint} shows  that expanding the network to avoid the congestion is financially prohibitive, which encourages us to develop stochastic control policies to gain gas network reliability and efficiency in a short run. 

Since the prediction of gas extractions involves errors, a gas network optimization problem has been addressed using the methods from robust optimization \cite{ben2009robust}, scenario-based and chance-constrained stochastic programming~\cite{shapiro2014lectures}. Besides forecasts, they require a network response model to uncertainty, i.e., the mapping from random forecast errors to the network state. The robust solutions \cite{Roald2020} optimize the network response to ensure the feasibility within robust uncertainty sets, but result in overly conservative operational costs. To alleviate the conservatism, scenario-based stochastic programs \cite{Zavala2014} optimize the network response to provide the minimum expected cost and ensure feasibility within a finite number of discrete scenarios. The major drawback of robust and scenario-based programs is their ignorance of the network state within the prescribed uncertainty set or outside the chosen scenarios. The chance-constrained programs \cite{ordoudis2018market}, \cite{Ratha2020}, in turn, yield an optimized network response across the entire forecast error distribution (or a family of those \cite{ordoudisenergy}), thus resulting in more advanced prediction and control of uncertain network state.

This work advocates the application of chance-constrained programming to the optimal natural gas network control under uncertainty. By optimal control, we imply the optimization of gas injection and pressure regulation policies that ensure gas flow feasibility and market efficiency for a given forecast error distribution. Towards this goal, we require a network response model with a strong analytic dependency between the network state and random forecast errors. Since natural gas flows are governed by non-convex equations, the design of network response models reduces to finding convex approximations. The work in~\cite[Chapter 6]{ordoudis2018market} enjoys the so-called controllable flow model~\cite{wang2017strategic}, which balances gas injection and uncertain extractions but disregards pressure variables. It thus does not permit policies for pressure control and corresponding financial remunerations. The work in~\cite{Ratha2020} preserves the integrity of system state variables and relies on the relaxation of non-convex equations. Although the relaxations are known to be tight~\cite{borraz2016convex,singh2020natural}, the results of \cite{Ratha2020} show that even a marginal relaxation gap yields a poor out-of-sample performance of the chance-constrained solution. Furthermore, the relaxations involve the integrality constraints to model bidirectional gas flows, which prevents extracting the dual solution and thus designing an optimal pricing scheme. One needs to introduce the unidirectional flow assumption to avoid integrality constraints, which is restrictive for gas networks under uncertainty \cite{Ratha2020}. 

This work bypasses the simplifying assumptions on network operations through the linearization of the non-convex natural gas equations, and provides a convex stochastic network optimization problem with performance guarantees. The problem ensures the real-time gas flow feasibility, enables the control of network state variability, and provides an efficient pricing scheme. Specifically, we make the following contributions:
\begin{enumerate}
    \item We propose stochastic control policies for gas injections and pressure regulation rates that provide real-time control inputs for network operators. Through linearization, we describe the uncertain state variables, such as nodal pressures and flow rates as affine functions of control inputs; thus capturing the dependency of the uncertain network state on operator's decisions. To establish performance guarantees, we provide a sample-based method to bound approximation errors induced due to linearization.
    \item We introduce a chance-constrained program to optimize the control policies and provide its computationally efficient second-order cone programming (SOCP) reformulation. The policy optimization ensures that the network state remains within network limits with a high probability and utilizes the statistical moments of the state variables to trade-off between the expected cost and the variance of the state variables.
    \item We propose a conic pricing scheme that remunerates network assets, i.e., gas suppliers, compressors and valves, for their contribution to uncertainty and variance control. Unlike the standard linear programming duality, the conic duality enables the decomposition of revenue streams associated with the coupling chance-constraints. We analyze the economic properties of the conic pricing scheme, e.g. revenue adequacy and cost recovery, and make them conditioned on the network design.
\end{enumerate}

At the operational planning stage, the optimized policies provide the best approximation (up to forecast quality) of the real-time control actions. They can be augmented into preoperational routines of network operators within the deterministic steady-state \cite{singh2020natural} or transient \cite{zlotnik2019optimal,zlotnik2019pipeline} gas models in the form of gas injection and pressure regulation set-points, while providing the strong foundation for necessary financial remunerations. We corroborate the effectiveness of the proposed policies using a $48$-node natural gas network. 

\textit{Outline:} Section \ref{sec:preliminaries} explains the gas network modeling, while Section \ref{sec:uncertainty} describes the stochastic network optimization, control policies and tractable reformulations. Section \ref{sec:pricing} introduces the pricing scheme and its theoretical properties. Section \ref{sec:experiments} provides numerical experiments, and Section \ref{sec:outlook} concludes. All proofs are relegated to Appendix. 

\textit{Notation:} Operation $\circ$ is the element-wise vector (matrix) product. Operator $\text{diag}[x]$ returns an $n\times n$ diagonal matrix with elements of vector $x\in\mathbb{R}^{n}$. For a $n\times n$ matrix $A$, $[A]_{i}$ returns an $i^{\text{th}}$ row $(1\times n)$ of matrix $A$, $\col{A}_{i}$ returns an $i^{\text{th}}$ column $(n\times 1)$ of matrix $A$, and $\text{Tr}[A]$ returns the trace of matrix $A$. Symbol $^{\top}$ stands for transposition, vector $\mathbb{1}$ $(\mathbb{0})$ is a vector of ones (zeros), and $\norm{\cdot}$ denotes the Euclidean norm.  

\section{Preliminaries} \label{sec:preliminaries}
\subsection{Gas Network Equations}
A natural gas network is modeled as a directed graph comprising a set of nodes $\mathcal{N}=\{1,\dots,N\}$ and a set of edges $\mathcal{E}=\{1,\dots,E\}$. Nodes represent the points of gas injection, extraction or network junction, while edges represent pipelines. Each edge is assigned a direction from sending node $n$ to receiving node $n'$, i.e., if $(n,n')\in\mathcal{E}$, then $(n',n)\notin\mathcal{E}$. The graph may contain cycles, while parallel edges and self-loops should not exist. The graph topology is described by a node-edge incidence matrix $A\in\mathbb{R}^{N\times E}$, such that
\begin{align*}
    A_{k\ell} = 
    \left\{
    \begin{array}{rl}
        +1, & \text{if}\;k=n  \\
        -1, & \text{if}\;k=n' \\
        0, & \text{otherwise}
    \end{array}
    \right.\quad\forall\ell=(n,n')\in\mathcal{E}.
\end{align*}
Let $\varphi\in\mathbb{R}^{E}$ be a vector of gas flow rates and let $\delta\in\mathbb{R}_{+}^{N}$ be a vector of gas extractions, which must be satisfied by the gas injections $\vartheta\in\mathbb{R}^{N}$ across the network given their injection limits $\underline{\vartheta},\overline{\vartheta}\in\mathbb{R}_{+}^{N}$. The gas conservation law is thus
\begin{align*}
    & A\varphi = \vartheta - \delta.
\end{align*}
The gas flow rates in network edges relate to the nodal pressures through non-linear, partial differential equations \cite{Thorley1987}. Under steady-state assumptions \cite{singh2020natural}, however, the flows are related to pressures through the Weymouth equation:
\begin{align*}
    \varphi_{\ell}|\varphi_{\ell}|=w_{\ell}\braceit{\varrho_{n}^2 - \varrho_{n'}^2},\quad\forall\ell=(n,n')\in\mathcal{E},
\end{align*}
where $\varrho\in\mathbb{R}^{N}$ is a vector of pressures contained within technical limits $\underline{\varrho},\overline{\varrho}\in\mathbb{R}_{+}^{N}$, and $w\in\mathbb{R}_{+}^{E}$ are constants that encode the friction coefficient and geometry of pipelines. To avoid non-linear pressure drops, let $\pi_{n}=\varrho_{n}^2$ be the squared pressure at node $n$ with limits $\underline{\pi}_{n}=\underline{\varrho}_{n}^2$ and $\overline{\pi}_{n}=\overline{\varrho}_{n}^2$. 

To support the desired nodal pressures, the gas network operator regulates the pressure using {\it active} pipelines $\mathcal{E}_{a}\subset\mathcal{E}$, which host either compressors $\mathcal{E}_{c}\subset\mathcal{E}_{a}$ or valves $\mathcal{E}_{v}\subset\mathcal{E}_{a}$, assuming $\mathcal{E}_{c}\cap\mathcal{E}_{v}=\emptyset$. These network assets respectively increase and decrease the gas pressure along their corresponding edges. To rewrite the gas conservation law and Weymouth equation  accounting for these components, let $\kappa\in\mathbb{R}^{E}$ be a vector of pressure regulation variables. Pressure regulation is non-negative $\kappa_{\ell}\geqslant0$ for every compressor edge $\ell\in\mathcal{E}_{c}$ and it is non-positive $\kappa_{\ell}\leqslant0$ for every valve edge $\ell\in\mathcal{E}_{v}$. This information is encoded in the pressure regulation limits $\underline{\kappa},\overline{\kappa}\in\mathbb{R}^{E}$. Pressure regulation involves an additional extraction of the gas mass to fuel active pipelines. Let matrix $B\in\mathbb{R}^{N\times E}$ relate the active pipelines to their sending nodes accounting for conversion factors, i.e., 
\begin{align*}
    B_{k\ell} = 
    \left\{
    \begin{array}{rl}
        b_{\ell}, & \text{if}\;k=n,\;k\in\mathcal{E}_{c}  \\
        -b_{\ell}, & \text{if}\;k=n,\;k\in\mathcal{E}_{v}  \\
        0, & \text{otherwise}
    \end{array}
    \right.\quad\forall\ell=(n,n')\in\mathcal{E},
\end{align*}
where $b_{\ell}$ is a conversion factor from the gas mass to the pressure regulation rate. The network equations become
\begin{subequations}
\begin{align}
     &A\varphi=\vartheta-B\kappa - \delta,\label{eq:conservation_low_vec}\\
     &\varphi\circ|\varphi|=\text{diag}[w](A^{\top}\pi + \kappa),\label{eq:weymouth_eq_vec}\\
     &\varphi_{\ell}\geqslant 0,\;\forall\ell\in\mathcal{E}_{a}.\label{eq:uni_direction}
\end{align}
\end{subequations}
Here, the gas extraction $B\kappa$ by compressor and valve edges in \eqref{eq:conservation_low_vec} is always non-negative. Equation \eqref{eq:weymouth_eq_vec} is the Weymouth equation in a vector form that accounts for both pressure loss and pressure regulation. The absolute value operator in \eqref{eq:weymouth_eq_vec} is understood element-wise. Finally, equality \eqref{eq:uni_direction} enforces the unidirectional condition for the gas flow in active pipelines, because they permit the gas flow only in one direction. 

\subsection{Deterministic Gas Network Optimization}
The gas network optimization seeks the minimum of gas injection costs while satisfying gas flow equations and network limits. Let $c_{1}\in\mathbb{R}_{+}^{N}$ and $c_{2}\in\mathbb{R}_{+}^{N}$ be the coefficients of a quadratic gas injection cost function. With a perfect extraction forecast, the deterministic gas network optimization is
\begin{subequations}\label{model:det}
\begin{align}
    \minimize{{\vartheta,\kappa,\varphi,\pi}}\quad&c_{1}^{\top}\vartheta + \vartheta^{\top}\text{diag}[c_{2}]\vartheta \label{det:obj}\\
    \st\quad&A\varphi=\vartheta-B\kappa - \delta, \label{det:conservation_low}\\
    &\varphi\circ|\varphi|=\text{diag}[w](A^{\top}\pi + \kappa),\label{det:weymouth_eq}\\
    &\underline{\pi}\leqslant\pi\leqslant\overline{\pi},\;\underline{\vartheta}\leqslant\vartheta\leqslant\overline{\vartheta},\label{det:lim_first}\\
    &\underline{\kappa}\leqslant\kappa\leqslant\overline{\kappa},\;\varphi_{\ell}\geqslant 0,\;\forall\ell\in\mathcal{E}_{a}.\label{det:lim_last}
\end{align}
\end{subequations}

Despite the non-convexity of \eqref{model:det}, it has been solved successfully using algorithmic solvers \cite{singh2019natural,singh2020natural} or general-purpose solvers \cite{wachter2006implementation} when all optimization parameters are known. These solvers no longer apply when the parameters are uncertain, because one needs to establish a convex dependency of optimization variables on uncertain parameters \cite{nemirovski2007convex}. This convex dependency is established in this work by means of the linearization of the Weymouth equation \eqref{det:weymouth_eq}. 

\subsection{Linearization of the Weymouth Equation}
Let $\mathcal{W}(\varphi,\pi,\kappa)=\mathbb{0}$ denote the non-convex constraint \eqref{det:weymouth_eq}, and let $\mathcal{J}(x)\in\mathbb{R}^{E\times n}$ denote the Jacobian of \eqref{det:weymouth_eq} w.r.t. an arbitrary vector $x\in\mathbb{R}^{n}$. The relation between the gas flow rates, nodal pressures, and pressure regulation rates can thus be approximated by the first-order Taylor series expansion:
\begin{align}
    \mathcal{W}(\varphi,\pi,\kappa)\approx& \mathcal{W}(\mathring{\varphi},\mathring{\pi},\mathring{\kappa}) 
    +\mathcal{J}(\mathring{\varphi})(\varphi-\mathring{\varphi}) \nonumber\\
    & + \mathcal{J}(\mathring{\pi})(\pi-\mathring{\pi})
    +\mathcal{J}(\mathring{\kappa})(\kappa-\mathring{\kappa})
    = \mathbb{0},\label{Taylor_expansion}
\end{align}
where $(\mathring{\varphi},\mathring{\pi},\mathring{\kappa})$ is a stationary point retrieved by solving non-convex problem \eqref{model:det}. As $\mathcal{W}(\mathring{\varphi},\mathring{\pi},\mathring{\kappa})=\mathbb{0}$ at a stationary point, equation \eqref{Taylor_expansion} implies the affine relation:
\begin{align}
    \varphi-\mathring{\varphi} &= 
    \mathcal{J}(\mathring{\varphi})^{-1}\mathcal{J}(\mathring{\pi})(\mathring{\pi}-\pi) 
    + \mathcal{J}(\mathring{\varphi})^{-1}\mathcal{J}(\mathring{\kappa})(\mathring{\kappa}-\kappa) \nonumber \\
    \Leftrightarrow
    \varphi &= 
    \underbrace{\mathcal{J}(\mathring{\varphi})^{-1}(\mathcal{J}(\mathring{\pi})\mathring{\pi} + \mathcal{J}(\mathring{\kappa})\mathring{\kappa}) + \mathring{\varphi}}_{\gamma_{1}(\mathring{\varphi},\mathring{\pi},\mathring{\kappa})} \nonumber \\
    &\quad\quad\quad\underbrace{-\mathcal{J}(\mathring{\varphi})^{-1}\mathcal{J}(\mathring{\pi})}_{\gamma_{2}(\mathring{\varphi},\mathring{\pi})}\pi
    \underbrace{-\mathcal{J}(\mathring{\varphi})^{-1}\mathcal{J}(\mathring{\kappa})}_{\gamma_{3}(\mathring{\varphi},\mathring{\kappa})}\kappa \nonumber\\
    \Leftrightarrow\varphi &= \gamma_{1}(\mathring{\varphi},\mathring{\pi},\mathring{\kappa}) + \gamma_{2}(\mathring{\varphi},\mathring{\pi})\pi + \gamma_{3}(\mathring{\varphi},\mathring{\kappa})\kappa, \label{flow_pressure_lin}
\end{align}    
where $\gamma_{1}\in\mathbb{R}^{E}$, $\gamma_{2} \in \mathbb{R}^{E \times N}$ and $\gamma_{3}\in\mathbb{R}^{E\times E}$ are coefficients encoding the sensitivity of gas flow rates to pressures and pressure regulation rates. These coefficients depend on the stationary point. For notational convenience, this dependency is dropped but always implied. In what follows, the Greek letter $\gamma$ denotes sensitivity coefficients and their transformations. 
\begin{remark}[Reference node]\label{rem:mul_ref_node}
    Since $\text{rank}(\gamma_{2}) = N-1$, system \eqref{flow_pressure_lin} is rank-deficient. Since the graph is connected, we have $E>N-1$, thus resulting in infinitely many solutions to system \eqref{flow_pressure_lin}. A unique solution is obtained by choosing a reference node $(\text{r})$ and fixing the reference pressure $\pi_{\text{r}} = \mathring{\pi}_{\text{r}}$. The reference node does not host a variable injection or extraction, nor should be a terminal node of active pipelines. In practice, this is a node with a large and constant gas injection. 
\end{remark}

\section{Gas Network Optimization under Uncertainty}\label{sec:uncertainty}
\subsection{Chance-Constrained Formulation}
At the operational planning stage, well ahead of the real-time operations, the unknown gas extractions are modeled as
\begin{align}\label{eq:rand_extr}
\tilde{\delta}(\xi) = \delta + \xi,
\end{align}
where $\delta \in \mathbb{R}^{N}$ is the mean value of the gas withdrawal rates and $\xi \in \mathbb{R}^{N}$ is a vector of zero-mean random forecast errors. Equation \eqref{eq:rand_extr} suffices to model disturbances in gas extractions without an explicit modeling of gas consumption by gas-fired power plants in adjacent electrical power grids. We assume that the forecast error distribution $\mathbb{P}_{\xi}$ of $\xi$ and covariance $\Sigma=\mathbb{E}[\xi\xi^{\top}]$ can be estimated from the historical observations of electrical loads and renewable power generation, that are known to obey Normal, Log-Normal and Weibull distributions \cite{hasan2019existing}. Though, more complex distributions may be envisaged for double-bounded stochastic processes of interest.

Regardless of the type and parameters of the uncertainty distribution, the chance-constrained counterpart of the deterministic gas network optimization in \eqref{model:det} writes as
\begin{subequations}\label{model:cc}
\begin{align}
    &\minimize{\tilde{\vartheta},\tilde{\kappa},\tilde{\varphi},\tilde{\pi}}\quad\mathbb{E}^{\mathbb{P}_{\xi}}[c_{1}^{\top}\tilde{\vartheta}(\xi) + \tilde{\vartheta}(\xi)^{\top}\text{diag}[c_{2}]\tilde{\vartheta}(\xi)] \label{cc:obj}\\
    &\quad\st\nonumber\\
    &\mathbb{P}_{\xi}\left[\begin{aligned}
    &A\tilde{\varphi}(\xi)=\tilde{\vartheta}(\xi)-B\tilde{\kappa}(\xi)-\tilde{\delta}(\xi),\\
    &\tilde{\varphi}(\xi) = \gamma_{1} + \gamma_{2}\tilde{\pi}(\xi) + \gamma_{3}\tilde{\kappa}(\xi),\\
    &\tilde{\pi}_{\text{r}}(\xi)=\mathring{\pi}_{\text{r}}
    \end{aligned}\right]\overset{\text{a.s.}}{=}1, \label{cc:con_as}
    \\
    &\mathbb{P}_{\xi}\left[
    \begin{aligned}
    &\underline{\pi}\leqslant\tilde{\pi}(\xi)\leqslant\overline{\pi},\;\underline{\vartheta}\leqslant\tilde{\vartheta}(\xi)\leqslant\overline{\vartheta},\\
    &\underline{\kappa}\leqslant\tilde{\kappa}(\xi)\leqslant\overline{\kappa},\;\tilde{\varphi}_{\ell}(\xi)\geqslant 0,\;\forall\ell\in\mathcal{E}_{a}
    \end{aligned}
	\right]\geqslant 1-\varepsilon,\label{cc:con_cc}
\end{align}
\end{subequations}
which optimizes stochastic network variables $\tilde{\vartheta},\tilde{\kappa},\tilde{\varphi}$ and $\tilde{\pi}$ to minimize the expected value of the cost function \eqref{cc:obj} subject to probabilistic constraints. The almost sure (a.s.) constraint \eqref{cc:con_as} requires the satisfaction of the gas conservation law and linearized Weymouth equation with probability 1, while the chance constraint \eqref{cc:con_cc} ensures that the real-time pressures together with the injection, pressure regulation and flow rates remain within their technical limits. The prescribed violation probability $\varepsilon\in(0,1)$ reflects the risk tolerance of the gas network operator towards the violation of network limits. 

\subsection{Control Policies and Network Response Model}
The chance-constrained problem \eqref{model:cc} is computationally intractable as it constitutes an infinite-dimensional optimization problem. To overcome its complexity, it has been proposed to approximate its solution by optimizing stochastic variables as affine, finite-dimensional functions of the random variable \cite{ben2009robust}. This functional dependency constitutes the model of the gas network response to uncertainty.

The \textit{explicit} dependency on uncertainty is enforced on the controllable variables through the following affine policies
\begin{subequations}\label{system_response}
\begin{align}
    \tilde{\vartheta}(\xi) = \vartheta + \alpha\xi, \quad\tilde{\kappa}(\xi) = \kappa + \beta\xi, \label{eq:explicit_policy}
\end{align}
where $\vartheta$ and $\kappa$ are the nominal (average) response, while $\alpha\in\mathbb{R}^{N \times N}$ and $\beta\in\mathbb{R}^{E \times N}$ are variable recourse decisions of the gas injections and pressure regulation by active pipelines, respectively. When optimized, policies \eqref{eq:explicit_policy} provide {\it control inputs} for the network operator to meet the realization of random forecast errors $\xi$. As the state variables, such as flow rates and pressures, are coupled with the controllable variables through stochastic equations \eqref{cc:con_as}, they \textit{implicitly} depend on uncertainty through the control inputs. 
\begin{lemma}\label{lem:1_state_variables_response}
Under control policies \eqref{eq:explicit_policy}, the random gas pressures and flow rates are given by affine functions
\begin{align}
    \tilde{\pi}(\xi) &= \pi + \breve{\gamma}_{2}(\alpha - \hat{\gamma}_{3}\beta -\text{diag}[\mathbb{1}])\xi,\label{eq:implicit_policy_pressure}\\
    \tilde{\varphi}(\xi) &= \varphi + (
    \grave{\gamma}_{2}(\alpha - \text{diag}[\mathbb{1}])
    - \grave{\gamma}_{3}\beta
    )\xi, \label{eq:implicit_policy_flow}
\end{align}
both including the nominal and random components, and where $\breve{\gamma}_{2},\hat{\gamma}_{2},\grave{\gamma}_{2},\hat{\gamma}_{3},\grave{\gamma}_{3}$ are constants of proper dimensions.
\end{lemma}
\end{subequations}
\medskip

Equations \eqref{system_response} constitute the desired model of the network response to uncertainty. The model is said to be admissible if the stochastic gas conservation law and linearized Weymouth equation in \eqref{cc:con_as} hold with probability 1, i.e., for any realization of random variable $\xi$. This is achieved as follows.  

\begin{lemma}\label{lem:2_recourse_eq}
The model of the gas network response \eqref{system_response} is admissible if the nominal and recourse variables obey
\begin{subequations}\label{eq:admisable_conditions}
\begin{align}
    &A\varphi=\vartheta - B\kappa - \delta \label{eq:admisable_conditions_1}\\
    &(\alpha - B\beta)^{\top}\mathbb{1} = \mathbb{1}, \label{eq:admisable_conditions_2}\\
    &\varphi = \gamma_{1} + \gamma_{2}\pi + \gamma_{3}\kappa, \label{eq:admisable_conditions_3}\\
    &\pi_{\text{r}} = \mathring{\pi}_{\text{r}},\;[\alpha]_{\text{r}}^{\top}=\mathbb{0},\;[\beta]_{\text{r}}^{\top}=\mathbb{0}. \label{eq:admisable_conditions_4}
\end{align}
\end{subequations}
\end{lemma}
\medskip

\begin{remark}\label{remark:distributuion_free}
    The model of the gas network response \eqref{system_response} does not make an assumption on the uncertainty distribution. 
\end{remark}

\subsection{Expected Cost Reformulation}

The expected value of the gas network cost function in \eqref{cc:obj} is computationally intractable as it involves an optimization of infinite-dimensional random variable $\tilde{\vartheta}(\xi)$. Under control policy \eqref{eq:explicit_policy}, however, we show that the computation of the expected cost reduces to solving an SOCP problem.

Due to definition of $\tilde{\vartheta}(\xi)$, function \eqref{cc:obj} rewrites as
\begin{align*}
    \mathbb{E}^{\mathbb{P}_{\xi}}[c_{1}^{\top}(\vartheta + \alpha\xi) + (\vartheta + \alpha\xi)^{\top}\text{diag}[c_{2}](\vartheta + \alpha\xi)],
\end{align*}
where the argument of the expectation operator is separable into nominal and random components. Due to the linearity of the expectation operator, it equivalently rewrites as
\begin{align*}
    c_{1}^{\top}\vartheta + \vartheta^{\top}\text{diag}[c_{2}]\vartheta +  \mathbb{E}^{\mathbb{P}_{\xi}}[c_{1}^{\top}\alpha\xi + (\alpha\xi)^{\top}\text{diag}[c_{2}]\alpha\xi].
\end{align*}
A zero-mean assumption made on distribution $\mathbb{P}_{\xi}$ factors out the first term under the expectation operator. The reformulation of the second term is made recalling that the expectation of the outer product of the zero-mean random variable yields its covariance, i.e., $\mathbb{E}[\xi\xi^{\top}]=\Sigma.$ Thus, the expected value of cost function \eqref{cc:obj} reduces to a computation of 
\begin{align*}
    c_{1}^{\top}\vartheta + \vartheta^{\top}\text{diag}[c_{2}]\vartheta + \text{Tr}[\alpha^{\top}\text{diag}[c_{2}]\alpha\Sigma],
\end{align*}
which is a convex quadratic function in variables $\vartheta$ and $\alpha$. To bring it to an SOCP form, let vectors $c^{\vartheta}\in\mathbb{R}^{N}$ and $c^{\alpha}\in\mathbb{R}^{N}$ substitute the quadratic terms of the gas injection and recourse costs. Moreover, let $F\in\mathbb{R}^{N\times N}$ be a factorization of covariance matrix $\Sigma$, such that $\Sigma=FF^{\top}$, and $\grave{c}_{2}\in\mathbb{R}^{N}$ be the factorization of vector $c_{2}$, such that $\text{diag}[c_{2}]=\grave{c}_{2}\grave{c}_{2}^{\top}$. Then, for any fixed values of nominal $\vartheta$ and recourse $\alpha$ decisions, the expected value of the cost is retrieved by solving the following SOCP problem
\begin{subequations}\label{model:cost_opt}
\begin{align}
    \minimize{c^{\vartheta},c^{\alpha}}\quad& c_{1}^{\top}\vartheta + \mathbb{1}^{\top}c^{\vartheta} + \mathbb{1}^{\top}c^{\alpha}\label{cost_opt_obj}\\
    \st\quad
    &\norm{\grave{c}_{2n}\vartheta_{n}}^{2}\leqslant c_{n}^{\vartheta},\;\forall n\in\mathcal{N}, \label{cost_opt_con_1}\\
    &\norm{F[\alpha]_{n}^{\top}c_{2n}}^{2}\leqslant c_{n}^{\alpha},\;\forall n\in\mathcal{N}, \label{cost_opt_con_2}
\end{align}
\end{subequations}
where \eqref{cost_opt_con_1} and \eqref{cost_opt_con_2} are rotated second-order cone constraints. Hence, the co-optimization of variables $\vartheta, \alpha, c^{\vartheta}$ and $c^{\alpha}$ results in the minimal expected cost. As problem \eqref{model:cost_opt} acts on a distribution-free response model (Remark \ref{remark:distributuion_free}), it does not require any assumption on the uncertainty distribution. 

\subsection{Variance of State Variables}

The optimization of response model \eqref{system_response} using the criterion of the minimum expected cost involves the risks of producing highly variable solutions for the state variables. See, for example, the evidences in the power system domain \cite{bienstock2019variance,mieth2020risk}. However, since the state variables \eqref{eq:implicit_policy_pressure} and \eqref{eq:implicit_policy_flow} are affine in control inputs, they can be optimized to provide the minimal-variance solution. To achieve the desired result, however, it is more suitable to optimize the standard deviations of the state variables as they admit conic formulations. 

Let $s^{\pi}\in\mathbb{R}^{N}$ and $s^{\varphi}\in\mathbb{R}^{E}$ be the variables modeling the standard deviations of pressures and flow rates, respectively. For any fixed values of recourse decisions $\alpha$ and $\beta$, the standard deviations of pressures and flows rates are retrieved by solving the following SOCP problem
\begin{subequations}\label{model:std_opt}
\begin{align}
    \minimize{s^{\pi},s^{\varphi}}\quad&\mathbb{1}^{\top}s^{\pi} + \mathbb{1}^{\top}s^{\varphi}\label{std_opt_obj}\\
    \st\quad&\norm{F[\breve{\gamma}_{2}(\alpha - \hat{\gamma}_{3}\beta -\text{diag}[\mathbb{1}])]_{n}^{\top}}\leqslant s_{n}^{\pi},\label{std_opt_con_1}\\
    &\norm{F[\grave{\gamma}_{2}(\alpha - \text{diag}[\mathbb{1}]) - \grave{\gamma}_{3}\beta]_{\ell}^{\top}}\leqslant s_{\ell}^{\varphi},\label{std_opt_con_2}\\
    &\forall n\in\mathcal{N}, \forall \ell\in\mathcal{E},\nonumber
\end{align}
\end{subequations}
where \eqref{std_opt_con_1} and \eqref{std_opt_con_2} are second-order cone constraints, which are tight at optimality. Therefore, the co-optimization of variables $\alpha, \beta, s^{\pi}$ and $s^{\varphi}$ yields the optimized system response \eqref{system_response} that ensures the minimal-variance solution for the state variables. We finally note that this co-optimization is also distribution-free. 

\subsection{Tractable Chance-Constrained Formulation}

It remains to reformulate the joint chance constraint \eqref{cc:con_cc} to attain a tractable reformulation. Given network response model \eqref{system_response}, one way to satisfy \eqref{cc:con_cc} is to enforce all its $N_{\mathsmaller{\leqslant}}$ inequalities on a finite number of samples from $\mathbb{P}_{\xi}$ \cite{campi2008exact}. The sample-based reformulation, however, does not explicitly parameterize the problem by the risk tolerance $\varepsilon$ of the network operator. We thus proceed by enforcing individual chance constraints with the explicit analytic parameterization of the risk tolerance through individual violation probabilities $\hat{\varepsilon}\in\mathbb{R}_{+}^{N_{\mathsmaller{\leqslant}}}$. This approach admits the Bonferroni approximation of the joint chance constraint in \eqref{cc:con_cc} when $\mathbb{1}^{\top}\hat{\varepsilon}\leqslant\varepsilon$. The joint feasibility guarantee is provided even when the choice of the individual violation probabilities is sub-optimal \cite{xie2019optimized}, e.g. $\hat{\varepsilon}_{i}=\frac{\varepsilon}{N_{\mathsmaller{\leqslant}}},\;\forall i=1,\dots, N_{\mathsmaller{\leqslant}}$. 

From \cite{nemirovski2007convex} we know that a scalar chance constraint 
\begin{subequations}\label{analytic_cc_reformulation}
\begin{align}
    \mathbb{P}_{\xi}[\xi^{\top}x\leqslant b]\geqslant 1-\hat{\varepsilon}\label{analytic_cc_reformulation_1}
\end{align}
analytically translates into the second-order cone constraint 
\begin{align}
    z_{\hat{\varepsilon}}\norm{Fx}\leqslant b - \mathbb{E}_{\xi}[\xi^{\top}x],\label{analytic_cc_reformulation_2}
\end{align}
\end{subequations}
where $z_{\hat{\varepsilon}}\geqslant0$ is a safety parameter in the sense of \cite{nemirovski2007convex}, and the left-hand side of \eqref{analytic_cc_reformulation_2} is the margin that ensures constraint feasibility given the parameters of the forecast errors distribution. Consequently, larger safety parameter $z_{\hat{\varepsilon}}$ improves system security. The choice of $z_{\hat{\varepsilon}}$ depends on the knowledge about distribution $\mathbb{P}_{\xi}$ \cite{nemirovski2007convex}, yet it always increases as the risk tolerance $\hat{\varepsilon}$ reduces. 

Given the  network response model \eqref{system_response} and the reformulations in \eqref{eq:admisable_conditions}--\eqref{analytic_cc_reformulation}, a computationally tractable version of stochastic problem \eqref{model:cc} with the variance awareness formulates as the following SOCP problem:  
\begin{subequations}\label{model:cc_ref}
\begin{align}
    \minimize{\mathcal{P}}\;\;c_{1}^{\top}&\vartheta + \mathbb{1}^{\top}c^{\vartheta} + \mathbb{1}^{\top}c^{\alpha} + \psi^{\pi\top}s^{\pi} + \psi^{\varphi\top}s^{\varphi}\label{cc_ref_obj}\\
    \st\;
    \lambda^{c}\colon
    &A\varphi=\vartheta - B\kappa - \delta, 
    \label{cc_ref_con_1}\\
    \lambda^{r}\colon
    &(\alpha - B\beta)^{\top}\mathbb{1} = \mathbb{1}, 
    \label{cc_ref_con_2}\\
    \lambda^{w}\colon
    &\varphi = \gamma_{1} + \gamma_{2}\pi + \gamma_{3}\kappa, \;\pi_{\text{r}} = \mathring{\pi}_{\text{r}},
    \label{cc_ref_con_3}\\
    \lambda_{n}^{\pi}\colon
    &\norm{F[\breve{\gamma}_{2}(\alpha - \hat{\gamma}_{3}\beta -\text{diag}[\mathbb{1}])]_{n}^{\top}}\leqslant s_{n}^{\pi},
    \label{cc_ref_con_4}\\
    \lambda_{\ell}^{\varphi}\colon
    &\norm{F[\grave{\gamma}_{2}(\alpha - \text{diag}[\mathbb{1}]) - \grave{\gamma}_{3}\beta]_{\ell}^{\top}}\leqslant s_{\ell}^{\varphi},
    \label{cc_ref_con_5}\\
    \lambda_{n}^{\overline{\pi}}\colon
    &z_{\hat{\varepsilon}}\norm{F[\breve{\gamma}_{2}(\alpha - \hat{\gamma}_{3}\beta -\text{diag}[\mathbb{1}])]_{n}^{\top}}\leqslant\overline{\pi}_{n}-\pi_{n},
    \label{cc_ref_con_6}\\
    \lambda_{n}^{\underline{\pi}}\colon
    &z_{\hat{\varepsilon}}\norm{F[\breve{\gamma}_{2}(\alpha - \hat{\gamma}_{3}\beta -\text{diag}[\mathbb{1}])]_{n}^{\top}}\leqslant\pi_{n}-\underline{\pi}_{n},
    \label{cc_ref_con_7}\\
    \lambda_{\ell}^{\underline{\varphi}}\colon
    &z_{\hat{\varepsilon}}\norm{F[\grave{\gamma}_{2}(\alpha - \text{diag}[\mathbb{1}]) - \grave{\gamma}_{3}\beta]_{\ell}^{\top}}\leqslant\varphi_{\ell},^{*}
    \label{cc_ref_con_8}\\
    &\textcolor{white}{z_{\hat{\varepsilon}}}\norm{\grave{c}_{2n}\vartheta_{n}}^{2}\leqslant c_{n}^{\vartheta},
    \label{cc_ref_con_9}\\
    &\textcolor{white}{z_{\hat{\varepsilon}}}\norm{F\grave{c}_{2n}[\alpha]_{n}^{\top}}^{2}\leqslant c_{n}^{\alpha},
    \label{cc_ref_con_10}\\
    &z_{\hat{\varepsilon}}\norm{F[\alpha]_{n}^{\top}}\leqslant\overline{\vartheta}_{n}-\vartheta_{n},
    \label{cc_ref_con_11}\\
    &z_{\hat{\varepsilon}}\norm{F[\alpha]_{n}^{\top}}\leqslant\vartheta_{n}-\underline{\vartheta}_{n},
    \label{cc_ref_con_12}\\
    &z_{\hat{\varepsilon}}\norm{F[\beta]_{\ell}^{\top}}\leqslant\overline{\kappa}_{\ell}-\kappa_{\ell},
    \label{cc_ref_con_13}\\
    &z_{\hat{\varepsilon}}\norm{F[\beta]_{\ell}^{\top}}\leqslant\kappa_{\ell}-\underline{\kappa}_{\ell},
    \label{cc_ref_con_14}\\
    &\forall n\in\mathcal{N},\;\forall \ell\in\mathcal{E},\;^{*}\forall \ell\in\mathcal{E}_{a},\nonumber
\end{align}
\end{subequations}
in variables $\mathcal{P}=\{\vartheta, \kappa, \varphi, \pi, \alpha, \beta, c^{\vartheta}, c^{\alpha}, s^{\pi}, s^{\varphi}\}$. Problem \eqref{model:cc_ref} optimizes the system response model \eqref{system_response} to meet a trade-off between the expected cost and the standard deviation of the state variables up to the given penalties $\psi^{\pi}\in\mathbb{R}_{+}^{N}$ and $\psi^{\varphi}\in\mathbb{R}_{+}^{E}$ for pressures and gas flow rates, respectively. Notice, that the constraints on the optimal recourse with respect to the reference node in \eqref{eq:admisable_conditions_4} are implicitly accounted for through the conic constraints on the gas injection and pressure regulation \eqref{cc_ref_con_11}--\eqref{cc_ref_con_14}. 

In formulation  \eqref{model:cc_ref}, the Greek letters $\lambda$ denote the dual variables of the coupling constraints. In the next Section \ref{sec:pricing}, we invoke the SOCP duality theory to establish an efficient pricing scheme for gas networks under uncertainty.

\subsection{Approximation Errors and Performance Guarantees}\label{subsec:errors}
Lemma \ref{lem:1_state_variables_response} hypothesizes the linear dependency of state variables on random forecast errors. Although the linear dependency enables a computationally tractable chance-constrained optimization in \eqref{model:cc_ref}, it also leads to approximation errors due to non-convex relation between pressures, flows, and uncertain gas extraction rates. To ensure that the optimization of control policies in \eqref{eq:explicit_policy} makes use of reliable state predictions, we develop \textit{a priori} worst-case performance guarantees that the approximation errors do not exceed a certain threshold. 

Since gas network congestions are mostly explained by pressure limits, we specifically focus on approximation errors associated with stochastic pressure variables. Let $\tilde{\pi}^{\star}(\xi)$ be the vector of the optimized stochastic pressures in \eqref{eq:implicit_policy_pressure}, i.e., 
\begin{align}
    \tilde{\pi}^{\star}(\xi) = \pi^{\star} + \breve{\gamma}_{2}(\alpha^{\star} - \hat{\gamma}_{3}\beta^{\star} -\text{diag}[\mathbb{1}])\xi,
\end{align}
which models the linear dependency on the optimal solution of problem \eqref{model:cc_ref}, denoted by $\star$, and random forecast error $\xi$.  

For some realization $\xi$, let $\pi^{\star}(\xi)$ be the actual pressure variables under control inputs from the optimized policies
\begin{align}
    \tilde{\vartheta}^{\star}(\xi) = \vartheta^{\star} + \alpha^{\star}\xi, \quad\tilde{\kappa}^{\star}(\xi) = \kappa^{\star} + \beta^{\star}\xi,\label{eq:opt_control_inputs}
\end{align}
where the optimal values are from the solution of problem \eqref{model:cc_ref}. Pressure variables $\pi^{\star}(\xi)$ can be then retrieved by projecting the optimized control inputs from \eqref{eq:opt_control_inputs} to the non-convex feasible region specific to realization $\xi$, i.e,  by solving
\begin{subequations}\label{model:projection_v1}
\begin{align}
    \pi^{\star}(\xi)\in\underset{\vartheta,\kappa,\varphi,\pi}{\text{argmin}}\quad&\norm{\tilde{\vartheta}^{\star}(\xi) - \vartheta} + \norm{\tilde{\kappa}^{\star}(\xi) - \kappa}\label{projection_v1_obj}\\
    \st\quad&A\varphi=\vartheta-B\kappa - (\delta + \xi),\\
    &\varphi\circ|\varphi|=\text{diag}[w](A^{\top}\pi + \kappa),\\
    &\text{Constraints}\;\eqref{det:lim_first}-\eqref{det:lim_last}.
\end{align}
\end{subequations}
For any node $n\in\mathcal{N}$, the stochastic pressure approximation error can be then pre-computed as an Euclidean distance
\begin{align}
    \Delta\pi_{n}(\xi) = \norm{\tilde{\pi}_{n}^{\star}(\xi) - \pi_{n}^{\star}(\xi)}
\end{align}
between the approximation $\tilde{\pi}_{n}^{\star}(\xi)$ and the actual pressure variable $\pi_{n}^{\star}(\xi)$ for some forecast error realization $\xi$. 

To provide the worst-case bound on the approximation error, we formulate the following optimization problem
\begin{subequations}\label{model:error_intractable}
\begin{align}
    \minimize{t}\quad&t\\
    \st\quad&\Delta\pi_{n}(\xi)-t\leqslant0,\quad\forall\xi\in\mathbb{P}_{\xi}\label{model_error_intractable_con}
\end{align}
\end{subequations}
in single variable $t$, which identifies that realization $\xi$ from $\mathbb{P}_{\xi}$, that results in the largest distance between the linear and non-convex stochastic pressure spaces. Observe, however, that constraint \eqref{model_error_intractable_con} is infinite as it requires infinitely many samples from $\mathbb{P}_{\xi}$. Using a sample-based approach from \cite{calafiore2005uncertain}, we provide the following finite counterpart of \eqref{model:error_intractable}
\begin{subequations}\label{model:error_tractable}
\begin{align}
    \minimize{t}\quad&t\\
    \st\quad&\Delta\pi_{n}(\widehat{\xi}_{s})-t\leqslant0,\quad\forall s = 1,\dots,S,\label{model_error_tractable_con}
\end{align}
\end{subequations}
where $\widehat{\xi}_{s}$ is a discrete sample from $\mathbb{P}_{\xi}$, and constraint \eqref{model_error_tractable_con} is enforced on a finite $S$ number samples (sample complexity), which is chosen to provide probabilistic performance guarantees with high confidence, as per the following Lemma. 
\begin{lemma}[Adapted from Corollary 1 in \cite{calafiore2005uncertain}]\label{lem:sample_comp}
For some $p\in[0,1]$ and $v\in[0,1]$, if sample complexity $S$ is such that 
\begin{align*}
S\geqslant\frac{1}{pv}-1,    
\end{align*}
then with probability $(1-p)$ and confidence level $(1-v)$, the pressure approximation error at node $n$ under the linear law in (7b) will not exceed the optimal solution $t^{\star}$ of problem \eqref{model:error_tractable}. 
\end{lemma}

\section{Pricing Gas Networks under Uncertainty}\label{sec:pricing}
From program \eqref{model:cc_ref}, we know that network assets participate in the satisfaction of the gas network equations \eqref{cc_ref_con_1}--\eqref{cc_ref_con_3}, in state variance reduction \eqref{cc_ref_con_4}--\eqref{cc_ref_con_5}, and in ensuring the feasibility of the state variables \eqref{cc_ref_con_6}--\eqref{cc_ref_con_8}. In this section, we establish a pricing scheme that remunerates network assets based on the combination of the classic linear programming duality \cite{kantorovich1960mathematical,samuelson1952spatial} and the SOCP duality \cite{boyd2004convex,mieth2020risk}. We refer the interested reader to Appendix \ref{app:dualization} 
for a brief overview on SOCP duality. For presentation clarity, however, we should stress that for each second-order cone constraint in \eqref{cc_ref_con_4}--\eqref{cc_ref_con_8} with a dual variable $\lambda\in\mathbb{R}^{1}$ there exists a vector of dual prices $u\in\mathbb{R}^{N}$, corresponding component-wise to random vector $\xi\in\mathbb{R}^{N}$, such that $\norm{u}\leqslant\lambda$. With a set of prices $\lambda,u_{1},\dots,u_{N}$, each conic coupling constraint becomes separable, thus enabling the revenue decomposition associated with constraints \eqref{cc_ref_con_4}--\eqref{cc_ref_con_8}. 

We first show that the primal and dual solutions of program \eqref{model:cc_ref} solve a \textit{partial} competitive equilibrium. This equilibrium consists of a price-setting problem that seeks the optimal prices associated with the coupling constraints \eqref{cc_ref_con_4}--\eqref{cc_ref_con_8},  a set of profit-maximizing problems of gas suppliers $n\in\mathcal{N}$, active pipelines $\ell\in\mathcal{E}_{a}$, and a rent-maximization problem solved by the network operator, as we establish in the proof of the following result; see Appendix \ref{proof_th_equilibrium} for details. Note, as program \eqref{model:cc_ref} does not model consumer preferences explicitly, we provide the results for partial equilibrium only.

\begin{table*}
\centering
\begin{minipage}{\textwidth}
\hrule
\medskip
\begin{subequations}\label{eq:revenues}
\begin{align}
    \mathcal{R}_{n}^{\text{sup}}&\triangleq
    \underbrace{
    \lambda_{n}^{c}\vartheta_{n} 
    }_{\text{\shortstack{nominal \\ balance}}}
    + 
    \underbrace{
    [\lambda^{r}]^{\top}[\alpha]_{n}^{\top}
    }_{\text{\shortstack{recourse \\ balance}}}
    + 
    \underbrace{
    z_{\hat{\varepsilon}}\braceit{\col{\breve{\gamma}_{2}}_{n}^{\top}(u^{\overline{\pi}} + u^{\underline{\pi}}) + \col{\grave{\gamma}_{2}}_{n}^{\top}u^{\underline{\varphi}}}F[\alpha]_{n}^{\top}
    }_{\text{\shortstack{gas pressure and flow limits}}}
    + 
    \underbrace{
    \braceit{\col{\breve{\gamma}_{2}}_{n}^{\top}u^{\pi} + \col{\grave{\gamma}_{2}}_{n}^{\top}u^{\varphi}}F[\alpha]_{n}^{\top}
    }_{\text{\shortstack{gas pressure and flow variance}}}
    \label{eq:revenues_inj}
    \\
    \mathcal{R}_{\ell}^{\text{act}}&\triangleq
    \underbrace{
    \braceit{\col{\gamma_{3}}_{\ell}^{\top}\lambda^{w}-\lambda^{c\top}\col{B}_{\ell}}\kappa_{\ell} 
    }_{\text{nominal pressure regulation}}
    - 
    \underbrace{
    \mathbb{1}^{\top}\col{B}_{\ell}\lambda^{r\top}[\beta]_{\ell}^{\top}
    }_{\text{\shortstack{recourse balance}}}
    - 
    \underbrace{
    z_{\hat{\varepsilon}}\braceit{\col{\breve{\gamma}_{2}\hat{\gamma}_{3}}_{\ell}^{\top} (u^{\overline{\pi}} + u^{\underline{\pi}}) + \col{\grave{\gamma}_{3}}_{\ell}^{\top} u^{\underline{\varphi}}} F[\beta]_{\ell}^{\top}
    }_{\text{\shortstack{gas pressure and flow limits}}}
    -
    \underbrace{
    \braceit{\col{\breve{\gamma}_{2}\hat{\gamma}_{3}}_{\ell}^{\top} u^{\pi} + \col{\grave{\gamma}_{3}}_{\ell}^{\top} u^{\varphi}} F[\beta]_{\ell}^{\top}
    }_{\text{\shortstack{gas pressure and flow variance}}}
    \label{eq:revenues_act}
    \\
    \mathcal{R}^{\text{rent}}&\triangleq
    \underbrace{
    \braceit{\lambda^{\underline{\varphi}\top} - \lambda^{w\top} - \lambda^{c\top} A}\varphi
    }_{\text{flow congestion rent}}
    +
    \underbrace{
    \braceit{\lambda^{w\top}\gamma_{2} + \lambda^{\underline{\pi}\top} -  \lambda^{\overline{\pi}\top}}\pi+\lambda^{\overline{\pi}\top}\overline{\pi} -\lambda^{\underline{\pi}\top}\underline{\pi}
    }_{\text{pressure congestion rent}}
    +
    \underbrace{
    \lambda^{\varphi\top} s^{\varphi} 
    +\lambda^{\pi\top}s^{\pi}
    }_{\text{variance rent}}
    \label{eq:revenues_rent}
    \\
    \mathcal{R}_{n}^{\text{con}}&\triangleq
    \underbrace{
    \lambda_{n}^{c}\delta_{n}
    }_{\text{\shortstack{nominal \\ balance}}}
    +
    \underbrace{
    \lambda_{n}^{r}
    }_{\text{\shortstack{recourse \\ balance}}}
    +
    \underbrace{
    z_{\hat{\varepsilon}}[F]_{n}\braceit{u^{\underline{\varphi}\top}\col{\grave{\gamma}_{2}}_{n} + (u^{\overline{\pi}}+ u^{\underline{\pi}})^{\top}\col{\breve{\gamma}_{2}}_{n}}
    }_{\text{\shortstack{gas pressure and flow limits}}}
    +
    \underbrace{
    [F]_{n}\braceit{u^{\varphi\top}\col{\grave{\gamma}_{2}}_{n} + u^{\pi\top} \col{\breve{\gamma}_{2}}_{n}}
    }_{\text{\shortstack{gas pressure and flow variance}}}
    \label{eq:revenues_con}
\end{align}
\end{subequations}
\medskip
\hrule
\end{minipage}
\vspace{-5mm}
\end{table*}

\begin{theorem}[Partial equilibrium payments]\label{th:equilibrium}
    Let $\mathcal{P}$ and $\mathcal{D}$ be the sets of the optimal primal and dual solutions of problem \eqref{model:cc_ref}, respectively. Then, both sets $\mathcal{P}$ and $\mathcal{D}$ solve a partial competitive network equilibrium with the following payments:
    \begin{itemize}
        \item Each gas supplier $n\in\mathcal{N}$ maximizes the expected profit when receiving the revenue of $\mathcal{R}_{n}^{\text{sup}}$ as in \eqref{eq:revenues_inj}. 
        \item Each active pipeline $\ell\in\mathcal{E}_{a}$ maximizes the expected profit when receiving the revenue of $\mathcal{R}_{\ell}^{\text{act}}$ as in \eqref{eq:revenues_act}. 
        \item The network operator minimizes the expected network congestion rent, which amounts to $\mathcal{R}^{\text{rent}}$ as in \eqref{eq:revenues_rent}. 
        \item The payment of each consumer $n\in\mathcal{N}$ is minimized when they are charged with $\mathcal{R}_{n}^{\text{con}}$  as in \eqref{eq:revenues_con}.
    \end{itemize}
\end{theorem}
\medskip

Similarly to a deterministic market settlement, the nominal gas injection or extraction is priced by associated locational marginal price $\lambda^{c}$, while the nominal pressure regulation is priced by the dual variable $\lambda^{w}$ of the Weymouth equation. The pricing scheme of Theorem \ref{th:equilibrium}, however, goes beyond the deterministic payments and provides three additional revenue streams for network assets \eqref{eq:revenues}. First, each network asset is paid with the dual variable $\lambda^{r}$ to remunerate its contribution to the feasibility of the gas network equations for any realization of uncertainty; see Lemma \ref{lem:2_recourse_eq}. The dual variables of the reformulated chance constraints \eqref{cc_ref_con_6}--\eqref{cc_ref_con_8} are used to compensate network assets for maintaining gas pressures and flow rates within network limits. Observe, this revenue stream is proportional to the safety parameter $z_{\hat{\varepsilon}}$, which increases as risk tolerance $\hat{\epsilon}$ reduces. The last revenue streams for network assets come from the satisfaction of the variance criteria set by the network operator. From the stationarity conditions \eqref{eq:stat_con_4_a} from Appendix \ref{proof_th_equilibrium}, the variance prices are $\lambda^{\pi}=\psi^{\pi}$ and $\lambda^{\varphi}=\psi^{\varphi}$, and from the SOCP dual feasibility condition \eqref{eq:dual_feasibility} from Appendix \ref{proof_th_equilibrium} we know that $\norm{[u^{\pi}]_{n}}\leqslant\lambda_{n}^{\pi},\;\norm{[u^{\varphi}]_{\ell}}\leqslant\lambda_{\ell}^{\varphi}$, $\forall n\in\mathcal{N}, \ell\in\mathcal{E}$. Thus, these revenue streams are proportional to the variance penalties $\psi^{\pi}$ and $\psi^{\varphi}$ set by the network operator. The consumer charges, motivated by their individual contributions to uncertainty and state variance, are explained similarly. Finally notice that, in contrast to the deterministic rent, revenue \eqref{eq:revenues_rent} additionally includes the variance control rent, which is non-zero whenever constraints \eqref{cc_ref_con_4}--\eqref{cc_ref_con_5} are binding, i.e., $\psi^{\pi},\psi^{\varphi}>0$.   

The results of Theorem \ref{th:equilibrium}, and thus the equivalence between the centralized optimization \eqref{model:cc_ref} and its equilibrium counterpart \eqref{model:profit_gas_injection}--\eqref{model:rent_minimization}, hold under certain assumptions. First, there exists at least one strictly feasible solution to SOCP problem \eqref{model:cc_ref} or to its dual counterpart to ensure that Slater's condition holds \cite{boyd2004convex}. Second, the market is perfectly competitive and the equilibrium agents act according to their true preferences, i.e., no exercise of market power. Finally, the information on the uncertainty distribution must be consistent among equilibrium problems \cite{dvorkin2019electricity}. Under these assumptions, we analyze the revenue adequacy and cost recovery of payments \eqref{eq:revenues} and make them conditioned on the network design.
\begin{corollary}[Revenue adequacy]\label{cor:revenue_adequcy}
    Let $\gamma_{1}=\mathbb{0}$ and $\underline{\pi}=\mathbb{0}$. Then, the payments established by Theorem \ref{th:equilibrium} are revenue adequate, i.e., $\mysum{n=1}{N}\mathcal{R}_{n}^{\text{con}}\geqslant\mysum{n=1}{N}\mathcal{R}_{n}^{\text{sup}} + \mysum{\ell=1}{E}\mathcal{R}_{\ell}^{\text{act}}. $
\end{corollary}
\medskip
As a result, the natural gas system does not incur a financial loss when the payments are distributed from consumers to network assets. The first condition in Corollary \ref{cor:revenue_adequcy} is motivated by the linearization of the Weymouth equation. If $\gamma_{1}\neq\mathbb{0}$, there exists an extra revenue term $\lambda^{w\top}\gamma_{1}$. As consumers are inelastic, this payment can be thus allocated to consumer charges, however its distribution among the customers remains an open question. Finally, the second condition in Corollary \ref{cor:revenue_adequcy} allows pressures to be zero at network nodes, which is too restrictive for practical purposes. In the next Section \ref{sec:experiments} we show that the revenue adequacy holds in practice even when this condition is not satisfied. 

The surplus of consumer payments in Corollary \ref{cor:revenue_adequcy} amounts to the congestion rent minimized by the network operator; see Appendix \ref{proof:rev_ad} for details. The consumer payments are thus implicitly minimized by problem \eqref{model:cc_ref} to only cover the congestion rent and compensate network assets for incurred costs. With our last result, we show that the cost recovery for network assets is also conditioned on the network design. 
\begin{corollary}[Cost recovery]\label{cor:cost_recovery}
     Let $\underline{\vartheta}=\mathbb{0},$ $\underline{\kappa}_{\ell}=0,\forall\ell\in\mathcal{E}_{c},$ and $\overline{\kappa}_{\ell}=0,\forall\ell\in\mathcal{E}_{v}$. Then, the payments of Theorem \ref{th:equilibrium} ensure cost recovery for suppliers and active pipelines, i.e., $\mathcal{R}_{\ell}^{\text{act}}\geqslant0,\forall\ell\in\mathcal{E}_{a},$ and $\mathcal{R}_{n}^{\text{sup}}-c_{1n}\vartheta_{n} - c_{n}^{\vartheta} - c_{n}^{\alpha}\geqslant0,\forall n\in\mathcal{N}.$
\end{corollary}
\medskip

\section{Numerical Experiments}\label{sec:experiments}

\begin{table*}
\caption{Deterministic versus Chance-Constrained Optimization of Control Policies}
\label{tab:netowrk_response}
\centering
\begin{tabular}{llcccccccc}
\toprule
\multirow{3}{*}{Parameter} & \multirow{3}{*}{Unit} & \multirow{3}{*}{\begin{tabular}[c]{@{}c@{}}Deterministic\\  control policies\end{tabular}} & \multicolumn{7}{c}{Chance-constrained control policies} \\
\cmidrule(lr){4-10}
 &  &  & \multirow{2}{*}{\begin{tabular}[c]{@{}c@{}}Variance-\\agnostic\end{tabular}} & \multicolumn{3}{c}{Pressure variance-aware, $\psi^{\pi}$} & \multicolumn{3}{c}{Flow variance-aware, $\psi^{\varphi}$} \\
\cmidrule(lr){5-7} \cmidrule(r){8-10}
 &  &  &  & $10^{-3}$ & $10^{-2}$ & $10^{-1}$ & $1$ & $10^{1}$ & $10^{2}$ \\
\midrule
Expected cost & \scalebox{0.85}{\$1000} & 80.9 & 82.5 (100\%) & 100.5\% & 105.6\% & 113.8\% & 100.1\% & 102.5\% & 112.6\% \\
$\sum_{n}\text{Var}[\tilde{\varrho}_{n}(\xi)]$ & \scalebox{0.85}{$\text{MPa}^2$} & 217.5 & 63.4 (100\%) & 44.2\% & 18.9\% & 12.8\% & 92.8\% & 46.7\% & 24.7\% \\
$\sum_{\ell}\text{Var}[\tilde{\varphi}_{\ell}(\xi)]$ & \scalebox{0.85}{$\text{BMSCFD}^2$} & 26.1 & 58.0 (100\%) & 83.4\% & 64.1\%  & 59.2\% & 93.4\% & 44.8\% & 25.9\% \\
\cmidrule(lr){1-10}
$\sum_{\ell \in \mathcal{E}_{c}}\sqrt{\kappa_{\ell}}$ & \scalebox{0.85}{$\text{kPa}$} & 1939 & 3914 & 3570 & 3734 & 3661 & 3914 & 4030 & 3888 \\
$\sum_{\ell \in \mathcal{E}_{v}}\sqrt{\kappa_{\ell}}$ & \scalebox{0.85}{$\text{kPa}$} & 0 & 0 & 0 & 150 & 576 & 0 & 1 & 500 \\
\cmidrule(lr){1-10}
Constraint inf. & \% & 53.7 & 0.04 & 0.02 & 0.02 & 0.02 & 0.03 & 0.02 & 0.03 \\
Average $\mathcal{P}_{\text{inj}}$ & \scalebox{0.85}{$\text{MMSCFD}$} & 960.91 & 0.01 & 0.03 & 0.02 & 0.02 & 0.02 & 0.04 & 0.04 \\
Average $\mathcal{P}_{\text{act}}$ & \scalebox{0.85}{$\text{kPa}$} & 121.68 & 0.19 & 0.08 & 0.10 & 0.05 &  0.28 & 0.04 & 0.04 \\
\bottomrule
\end{tabular}
\end{table*}

\begin{figure*}[t]
\centering
\includegraphics[width=\textwidth]{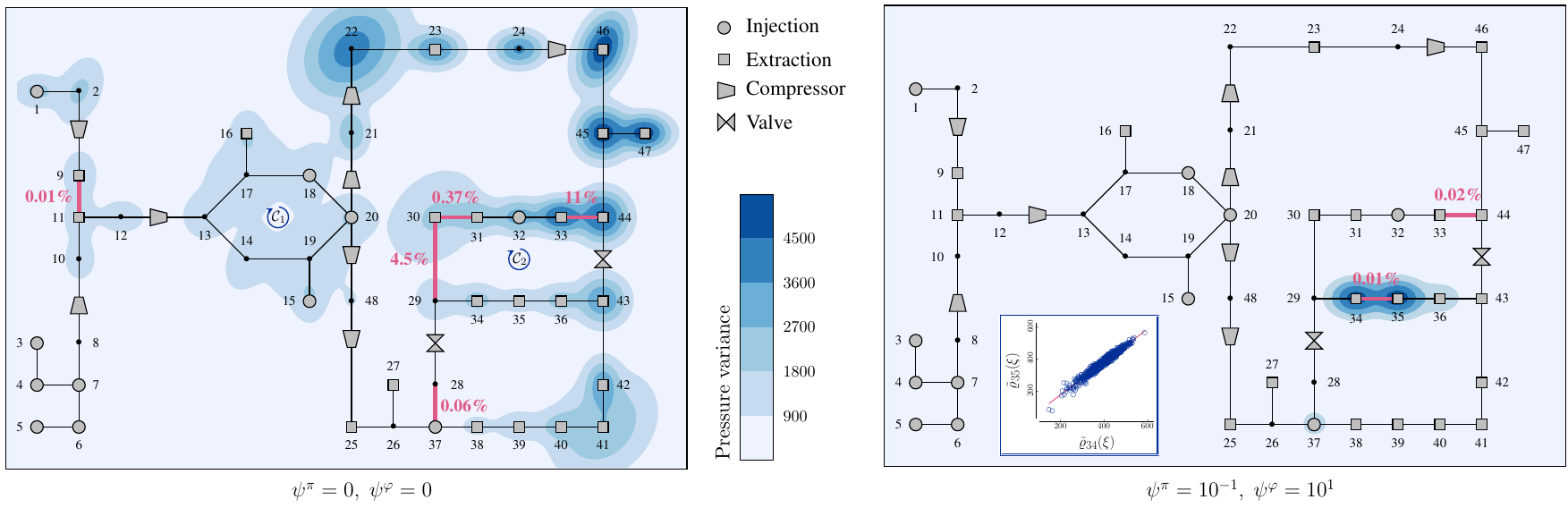}
\vspace{-5mm}
\caption{Comparison of the variance-agnostic (left) and the variance-aware (right) chance-constrained control policies in terms of the state variables variance for $\varepsilon=10\%$. The red values show the probability of flow reversal. The inset plot shows the correlation between the pressures at nodes $34$ and $35$.}
\label{fig:var_pen}
\vspace{-5mm}
\end{figure*}
\begin{figure}
\begin{center}
\resizebox{0.45\textwidth}{!}{%
\includegraphics[width=0.5\textwidth]{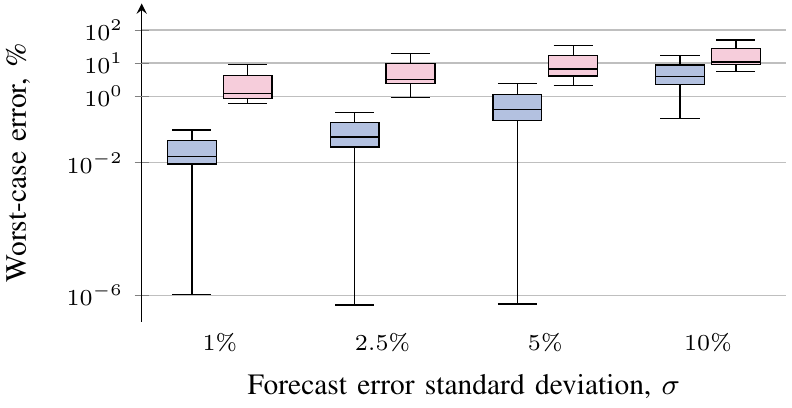}}
\end{center}
\vspace{-5mm}
\caption{The worst-case stochastic pressure approximation errors summarized across 48 nodes, for $p=v=0.9$ in Lemma \ref{lem:sample_comp}, and for different uncertainty penetration levels. The blue boxplots \includegraphics[width=0.25cm]{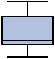} summarize error statistics when the chance-constrained control policies are implemented, and red boxplots \includegraphics[width=0.25cm]{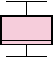} provide the summary of the deterministic solution.}
\label{fig:lin_error}
\end{figure}
\begin{figure}
\center
\includegraphics[width=0.45\textwidth]{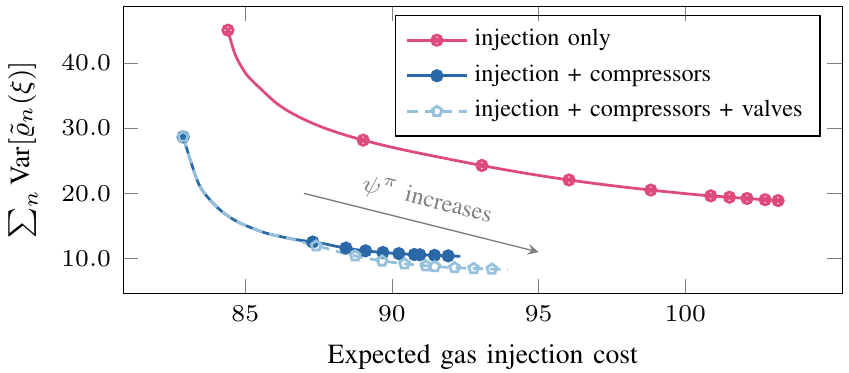}
\caption{Expected cost versus pressure variance for different assignments of control polices to network assets. Pressure penalty $\psi^{\pi}\in[10^{-3},10^{-1}]$.}
\label{fig:variance_cost_topology}
\vspace{-5mm}
\end{figure}
\begin{figure}
\center
\resizebox{0.45\textwidth}{!}{%
\includegraphics[width=0.5\textwidth]{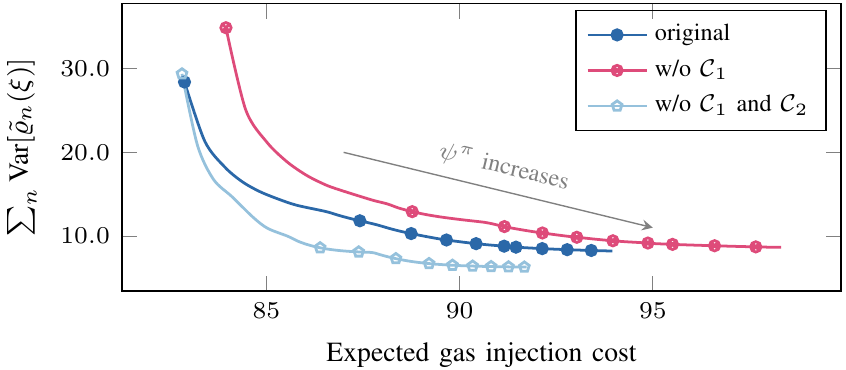}
}
\caption{Expected cost versus pressure variance under three network structures. Pressure penalty $\psi^{\pi}\in[10^{-3},10^{-1}]$.}
\label{fig:variance_cost_network_structure}
\end{figure}

We run numerical experiments using a $48$-node natural gas network depicted in Fig. \ref{fig:var_pen}. The network parameters are sourced from \cite{Wu2000} with a few modifications to enable large uncertainty and variability of gas extraction rates and provide more variance control opportunities. Specifically, the pressure limits at nodes $1$ and $3$ are homogenized with those at the rest of network nodes, two injections are added in the demand area at nodes $32$ and $37$, and two valves are installed in the pipelines connecting nodes $(28,29)$ and $(43,44)$. The $22$ gas extractions are modeled as $\tilde{\delta}(\xi) = \delta + \xi $, where $\delta$ is the nominal extraction rate reported in \cite{Wu2000} and $\xi$ is the zero-mean normally distributed forecast error. The safety parameter $z_{\hat{\varepsilon}}$ is thus the inverse CDF of the standard Gaussian distribution at $(1-\hat{\varepsilon})-$quantile \cite{nemirovski2007convex}. The standard deviation of each gas extraction is set to $10\%$ of the nominal rate. The joint constraint violation probability $\varepsilon$ is set to 1\% by default. To retrieve the stationary point in \eqref{flow_pressure_lin}, the non-convex problem \eqref{model:det} is solved for the nominal gas extraction rates using the Ipopt solver \cite{wachter2006implementation}. The repository \cite{companion2020dvorkin} contains the input data and code implementation in the JuMP package for Julia \cite{dunning2017jump}.  

\subsection{Analysis of the Optimized Network Response}
We first study the optimized gas network response to uncertainty under deterministic and chance-constrained control policies \eqref{eq:explicit_policy}. The deterministic policies are optimized by setting the safety factor $z_{\hat{\varepsilon}}$ in problem \eqref{model:cc_ref} to zero. The policies are compared in terms of the expected cost \eqref{cost_opt_obj}, the aggregated variance of gas pressures and flow rates  \scalebox{0.9}{$\sum_{n}\text{Var}[\tilde{\varrho}_{n}(\xi)]$} and \scalebox{0.9}{$\sum_{\ell}\text{Var}[\tilde{\varphi}_{\ell}(\xi)]$}, respectively, and the total pressure regulation by compressors \scalebox{0.9}{$\sum_{\ell \in \mathcal{E}_{c}}\sqrt{\kappa_{\ell}}$} and valves \scalebox{0.9}{$\sum_{\ell \in \mathcal{E}_{v}}\sqrt{\kappa_{\ell}}$}. Note, we discuss the natural pressure quantities, not their squared counterparts used in optimization. The policies are also compared in terms of network constraints satisfaction. We first sample control inputs from \eqref{system_response} for $S=1,000$ realizations of forecast errors and count the violations of network limits \eqref{cc:con_cc}. Second, we assess the quality of the control inputs \eqref{eq:explicit_policy} for the non-convex gas equations, by solving the projection problem \eqref{model:projection_v1} for all realizations $\xi_{s}, \forall s=1,\dots,S$. A control input is considered feasible if \eqref{projection_v1_obj} is zero for a given realization. To characterize this infeasibility numerically, consider the average metrics $\mathcal{P}_{\text{inj}}=\sum_{s}\norm{\tilde{\vartheta}(\xi_{s}) - \vartheta_{s}}/S$ for gas injections and $\mathcal{P}_{\text{act}}=\sum_{s}\norm{\tilde{\kappa}(\xi_{s}) - \kappa_{s}}/S$ for active pipelines.

The results are reported in Table \ref{tab:netowrk_response}. Disregarding uncertainty, the deterministic policies optimize the network operation for the nominal gas extraction rates and thus result in the minimum of cost at the operational planning stage. However, the produced control inputs are infeasible for most of the forecast error realizations. The projections $\mathcal{P}_{\text{inj}}$ and $\mathcal{P}_{\text{act}}$ of deterministic policies require the real-time correction of gas injections by $31.3$\% and the real-time correction of pressure regulation by active pipelines by $12.7$\% of the nominal rates on average. The chance-constrained policies, on the other hand, produce the control inputs that remain feasible with a probability at least $1-\varepsilon=99\%$ and require a minimal effort to restore the real-time gas flow feasibility. This real-time effort is non-zero due to approximation errors induced by linear pressure and flow equations of Lemma \ref{lem:1_state_variables_response}. Figure \ref{fig:lin_error} illustrates the \textit{worst-case} stochastic pressure approximation errors obtained according to the approach in Section \ref{subsec:errors}. The errors significantly depend on the amount of uncertainty: with probability $90\%$ and at high confidence, the errors approach 0\% for a small uncertainty penetration level ($\sigma=1\%$), and they will not exceed 5.8\% on average for the extremely large uncertainty penetration ($\sigma=10\%$). The errors under the deterministic solution, which ignores gas extraction uncertainty, are larger by at least an order of magnitude on average.  

Table \ref{tab:netowrk_response} further demonstrates that the variance-agnostic policy requires only a slight increase of the expected cost relative to the deterministic solution by $1.6$\%, while the variance-aware policies allow to trade-off the expected operational cost for the smaller variations of pressures and flow rates. The variance of gas pressures and flow rates can be reduced by $63.8$\% and $7.2$\%, respectively, without any substantial impact on the expected cost. Observe that the subsequent variance reduction is achieved also due to the activation of valves in two active pipelines, that are not operating in the deterministic and variance-agnostic solutions.  

With the density plots in Fig. \ref{fig:var_pen}, we demonstrate the uncertainty propagation through the network. The variance-agnostic solution results in the large pressure variance in the eastern part of the network with a large concentration of stochastic gas extractions. This solution further allows the probability of the gas flows reversal up to $11$\% for certain pipelines, thus making the prediction of flow directions difficult. The variance-aware solution with the joint penalization of pressures and flows variance, in turn, drastically reduces the variation of the state variables and localizes the most of the variation only at nodes $34$ and $35$. Failure to minimize the pressure variance at these two nodes is due to relatively large approximation errors compared to the rest of the nodes (see the top quantiles of blue boxplots in Fig. \ref{fig:lin_error}).
Although this variation remains large, the pressures at these nodes are highly correlated. Thus, by Weymouth equation \eqref{det:weymouth_eq}, the flow variance and the probability of flow reversal in edge $(34,35)$ remain small. 

Next, we analyze the contribution of network assets to the variance control through the cost-variance trade-offs in Fig. \ref{fig:variance_cost_topology}. The figure illustrates these trade-offs when the control policies are assigned to gas injections only ($\alpha\in\text{free},\beta=\mathbb{0}$), to gas injections and compressors ($\alpha,\beta\in\text{free},[\beta]_{\ell}^{\top}=\mathbb{0},\forall\ell\in\mathcal{E}_{v}$), and to all network assets including valves ($\alpha,\beta\in\text{free}$). Observe that the variance reduction is achieved more rapidly and at lower costs as more active pipelines are involved into uncertainty and variance control. Hence, the stochastic control becomes more available as the network operator deploys more pressure regulation action by compressors and valves.

Last, we analyze structural network impacts on the cost-variance trade-offs. We gradually break cycles $\mathcal{C}_{1}$ (by removing edges $(13,14)$ and $(14,19)$) and $\mathcal{C}_{2}$ (by removing edge $(29,30)$) in Fig. \ref{fig:var_pen} to change the network to a tree-like topology, leaving only those cycles that are mandatory for feasible operation. Figure \ref{fig:variance_cost_network_structure} summarizes the cost-variance trade-offs and points on the ambiguous role of network cycles. Breaking cycle $\mathcal{C}_{1}$ in the supply concentration area causes congestion, which prevents deploying western suppliers to minimize the pressure variance in the east, substantially increasing the cost of operations. On the other hand, the subsequent removal of cycle $\mathcal{C}_{2}$ weakens the graph connectivity and allows for more economical and more drastic variance reduction in the east. This agrees with equation \eqref{eq:implicit_policy_pressure}, which relates pressures and forecast errors through parameters $\breve{\gamma}_{2}$ and $\hat{\gamma}_{3}$, that encode graph connectivity. We notice, the trade-offs in Fig.  \ref{fig:variance_cost_topology} and \ref{fig:variance_cost_network_structure} motivate the problem of the variance-aware network design.

\subsection{Revenue Analysis}
Figure \ref{fig:revenue} depicts the total revenues of active pipelines and gas injections as well as the total charges of gas consumers. It further shows their decomposition into revenue streams defined by the pricing scheme in \eqref{eq:revenues}. Relative to the deterministic payments, the chance-constrained policies lead to a substantial increase in payments that further increase due to the variance awareness. Besides the nominal supply revenues, the chance-constrained policies produce the compensations for the uncertainty and variance control that together exceed deterministic payments by $37.3\%$. Moreover, the payments for the nominal supply under stochastic policies also increase due to several reasons. First, as shown in Table \ref{tab:netowrk_response}, the stochastic policies require a larger deployment of gas compressors and valves that extract an additional gas mass for fuel purposes, up to 4.2\% of the network demand, thus increasing the marginal cost of gas suppliers. Second, to provide the security margins for chance constraints \eqref{cc_ref_con_6}--\eqref{cc_ref_con_8} and \eqref{cc_ref_con_11}--\eqref{cc_ref_con_14}, the optimized policies require withholding less expensive injections from the purposes of the nominal supply. 
Last, with increasing assignments of penalty factors $\psi^{\pi}$ and $\psi^{\varphi}$, 
the optimality of the nominal injection cost is altered in the interest of reduced variance of state variables.  
Finally, the mismatch between the consumer charges and the revenues of gas injections and active pipelines is non-negative, thus satisfying the revenue adequacy in all three instances.  

\begin{figure}
\centering
\includegraphics[width=0.5\textwidth]{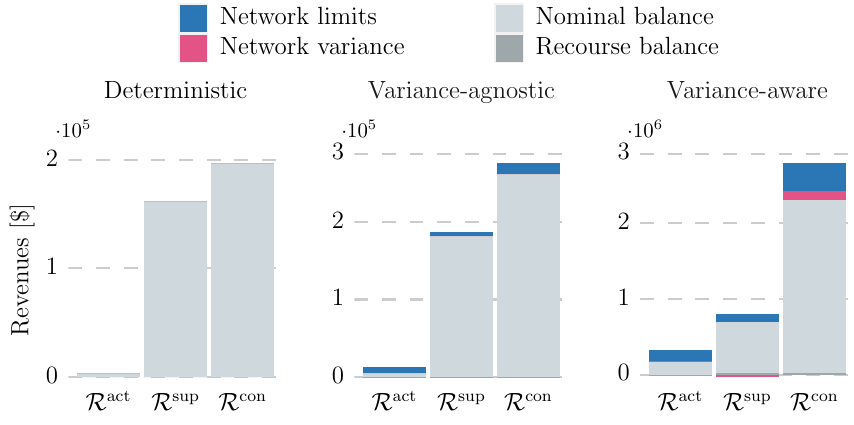}
\vspace{-5mm}
\caption{Total payments for active pipelines $\mathcal{R}^{\text{act}}$, suppliers $\mathcal{R}^{\text{sup}}$ and consumers $\mathcal{R}^{\text{con}}$ under deterministic, chance-constrained variance-agnostic and chance-constrained variance-aware ($\psi^{\pi}=0.1,\psi^\varphi=100$) policies.}
\label{fig:revenue}
\vspace{-5mm}
\end{figure}

\section{Conclusions \& Outlook}\label{sec:outlook}
This work has established the stochastic control policies and pricing scheme for the non-convex steady-state gas network operations under gas extraction uncertainty. The work offers an uncertainty- and variance-aware policy optimization that ensures the gas flow feasibility with a high probability and minimal variance of the state variables. Moreover, the work challenged the deterministic market settlement and offered financial remunerations to network assets for their contributions to uncertainty and variance control. 

The definition and optimization of gas storage control policies under uncertainty constitute the relevant direction for a future work. In addition, a price-responsive modeling of stochastic gas extraction rates, by means of co-optimization of control policies and the stochastic gas consumption models, is a valid research direction. This would lead, for example, to uncertainty- and variance-aware coordination and financial contracts between the gas and power network operators.

\appendix 
\subsection{Proof of Lemma \ref{lem:1_state_variables_response}}
The substitution of the linearized Weymouth equation from \eqref{cc:con_as} and policies \eqref{eq:explicit_policy} into the gas conservation law in \eqref{cc:con_as} yields stochastic pressures as
\begin{align*}
    &A\tilde{\varphi}(\xi)=\tilde{\vartheta}(\xi) - B\tilde{\kappa}(\xi) - \tilde{\delta}(\xi)\\
    &\Leftrightarrow
    A(\gamma_{1} + \gamma_{2}\tilde{\pi}(\xi) + \gamma_{3}(\kappa + \beta\xi)) \\
    &\qquad \qquad \qquad \qquad = \vartheta + \alpha\xi - B(\kappa + \beta\xi) - \delta - \xi\\
    &\Leftrightarrow 
    \underbrace{A\gamma_{2}}_{\hat{\gamma}_{2}}\tilde{\pi}(\xi) = 
    \underbrace{
    \vartheta - (B + A\gamma_{3})\kappa - \delta - A\gamma_{1}
    }_{\text{from}\;\eqref{det:conservation_low},\eqref{flow_pressure_lin}\;\colon A\gamma_{2}\pi = \hat{\gamma}_{2}\pi} \\
    &\quad\quad\quad\quad\quad\quad+ (\alpha - \underbrace{(B + A\gamma_{3})}_{\hat{\gamma}_{3}}\beta - \text{diag}[\mathbb{1}])\xi\\
    &\Leftrightarrow 
    \hat{\gamma}_{2}\tilde{\pi}(\xi) = \hat{\gamma}_{2}\pi + (\alpha - \hat{\gamma}_{3}\beta - \text{diag}[\mathbb{1}])\xi\\
    &\Leftrightarrow 
    \tilde{\pi}(\xi) = \pi + \hat{\gamma}_{2}^{-1}(\alpha - \hat{\gamma}_{3}\beta -\text{diag}[\mathbb{1}])\xi,
\end{align*}
where $\hat{\gamma}_{2} \in \mathbb{R}^{N \times N}$ and $\hat{\gamma}_{3}\in \mathbb{R}^{N \times E}$ are auxiliary constants. As $\hat{\gamma}_{2}=A\gamma_{2}$, it is only invertible for the tree network topology. For generality, consider a reference node (r), see Remark \ref{rem:mul_ref_node}, and let $\hat{\gamma}_{2\backslash\text{r}}$ be a reduced matrix $\hat{\gamma}_{2}$ without the $\text{r}^{\text{th}}$ row and column in $\hat{\gamma}_{2}$. The invertible counterpart of $\hat{\gamma}_{2}$ is
 \begin{align*}
    \breve{\gamma}_{2} = \begin{bmatrix}\hat{\gamma}_{2\backslash\text{r}}^{-1}&\mathbb{0}\\\mathbb{0}^{\top}&0\end{bmatrix},
 \end{align*}
and the stochastic pressures become
\begin{subequations}
\begin{align}
&\tilde{\pi}(\xi) = \pi + \breve{\gamma}_{2}(\alpha - \hat{\gamma}_{3}\beta -\text{diag}[\mathbb{1}])\xi,\label{proof_press_stoch}\\ 
&\pi_{\text{r}} = \mathring{\pi}_{\text{r}},\;[\alpha]_{\text{r}}^{\top}=\mathbb{0},\;[\beta]_{\text{r}}^{\top}=\mathbb{0},\label{proof_press_ref}
\end{align}
\end{subequations}
for an arbitrary network topology. Here, equation \eqref{proof_press_ref} is enforced to satisfy the reference node definition. 

To obtain the stochastic flow rates, substitute \eqref{proof_press_stoch} into the linearized Weymouth equation in \eqref{cc:con_as} and rearrange, i.e.,
\begin{align*}
    &\tilde{\varphi}(\xi) = \gamma_{1} + \gamma_{2}\tilde{\pi}(\xi) + \gamma_{3}\tilde{\kappa}(\xi) \nonumber\\
    &\Leftrightarrow
    \tilde{\varphi}(\xi) = 
    \underbrace{\gamma_{1} + \gamma_{2}\pi + \gamma_{3}\kappa}_{\text{from}\;\eqref{flow_pressure_lin}\;\colon\varphi} 
    + \underbrace{\gamma_{2}\breve{\gamma}_{2}}_{\grave{\gamma}_{2}}(\alpha - \text{diag}[\mathbb{1}])\xi \\
    &\quad\quad\quad\quad\quad- \underbrace{(\gamma_{2}\breve{\gamma}_{2}\hat{\gamma}_{3} - \gamma_{3})}_{\grave{\gamma}_{3}}\beta\xi\\
    &\Leftrightarrow
    \tilde{\varphi}(\xi) = \varphi + (\grave{\gamma}_{2}(\alpha - \text{diag}[\mathbb{1}]) + \grave{\gamma}_{3}\beta)\xi,
\end{align*}
where $\grave{\gamma}_{2}\in\mathbb{R}^{E\times N}$ and $\grave{\gamma}_{3} \in \mathbb{R}^{E \times E}$ are constants.

\subsection{Proof of Lemma \ref{lem:2_recourse_eq}}
Consider the stochastic gas conservation law in \eqref{cc:con_as}:
\begin{align*}
    A\tilde{\varphi}(\xi)=\tilde{\vartheta}(\xi)-B\tilde{\kappa}(\xi)-\tilde{\delta}(\xi).
\end{align*}
From the properties of the edge-node incidence matrix $A$, we know that $\mathbb{1}^{\top}A\tilde{\varphi}(\xi) = 0.$ By summing up $N$ equations above and by substituting equations \eqref{eq:explicit_policy}, we arrive to equation
\begin{align*}
    \mathbb{1}^{\top}\vartheta - \mathbb{1}^{\top}B\kappa - \mathbb{1}^{\top}\delta + \mathbb{1}^{\top}\alpha\xi - \mathbb{1}^{\top}B\beta\xi - \mathbb{1}^{\top}\xi = 0,
\end{align*}
which is separable into nominal and random components:
\begin{subequations}
\begin{align}
    &\mathbb{1}^{\top}\vartheta - \mathbb{1}^{\top}B\kappa - \mathbb{1}^{\top}\delta = 0,\label{app_aux_1}\\
    &\mathbb{1}^{\top}\alpha\xi - \mathbb{1}^{\top}B\beta\xi - \mathbb{1}^{\top}\xi = 0,\label{app_aux_2}
\end{align}
\end{subequations}
where equation \eqref{app_aux_1} is the deterministic gas conservation law, which is alternatively expressed through \eqref{eq:conservation_low_vec}, thus providing the first condition in \eqref{eq:admisable_conditions_1}. The second condition in \eqref{eq:admisable_conditions_2} is provided from \eqref{app_aux_2}, which holds for any realization of $\xi$ if the recourse variables $\alpha$ and $\beta$ obey $(\alpha - B\beta)^{\top}\mathbb{1} = \mathbb{1}$.

To obtain condition \eqref{eq:admisable_conditions_3}, substitute \eqref{system_response} into the stochastic linearized Weymouth equation in \eqref{cc:con_as}:
\begin{align*}
    \varphi =&\; \gamma_{1} + \gamma_{2}\pi + \gamma_{3}\kappa 
    - \alpha(\grave{\gamma}_{2} - \gamma_{2}\breve{\gamma}_{2})\xi \\
    &+ \beta(\grave{\gamma}_{3} - \gamma_{2}\breve{\gamma}_{2}\hat{\gamma}_{3} + \gamma_{3})\xi
    + (\grave{\gamma}_{2} - \gamma_{2}\breve{\gamma}_{2})\text{diag}[\mathbb{1}]\xi\\
    =&\; \gamma_{1} + \gamma_{2}\pi + \gamma_{3}\kappa,
\end{align*}
yielding a deterministic equation due to the definition of constants $\grave{\gamma}_{2}$ and $\grave{\gamma}_{3}$. Finally, the stochastic equation for the reference node is satisfied by equations \eqref{proof_press_ref}. 

\subsection{Dualization of Conic Constraints}\label{app:dualization}

\begin{subequations}
The results presented in this section are due to \cite[Chapter 5]{boyd2004convex}. Consider the SOCP problem of the form
\begin{align}
    \minimize{x\in\mathbb{R}^{n}}\quad c^{\top}x,\;
    \st\quad\norm{A_{i}x}\leqslant b_{i}^{\top}x,\quad\forall i=1,\dots,m,
\end{align}
with $c\in\mathbb{R}^{n}$, $A_{i}\in\mathbb{R}^{n_{i}\times n}$, $b_{i}\in\mathbb{R}^{n}$. 
To dualize the second-order cone constraint, we use the fact that for any pair $\lambda_{i}\in\mathbb{R}^{1}$ and $u_{i}\in\mathbb{R}^{n_{i}}$ it holds that
\begin{align}
\maximize{u_{i},\lambda_{i}:\atop\norm{u_{i}}\leqslant\lambda_{i}}&\;-u_{i}^{\top}A_{i}x-\lambda_{i}b_{i}^{\top}x 
= 
\maximize{\lambda_{i}\geqslant0}\;-\lambda_{i}(\norm{A_{i}x}- b_{i}^{\top}x) \nonumber \\
&= 
\left\{
\begin{matrix*}[l]
0, & \text{if}\;\norm{A_{i}x}\leqslant b_{i}^{\top}x,\\
-\infty, &\text{otherwise}.
\end{matrix*}
\right.
\label{eq:dualization_soc}
\end{align}
Therefore, the Lagrangian of the SOCP problem writes in variables $x\in\mathbb{R}^{n}$, $\lambda\in\mathbb{R}^{m}$ and $u\in\mathbb{R}^{n_{i}\times n}$ as 
\begin{align}
    \maximize{\norm{u_{i}}\leqslant\lambda_{i}}\quad\minimize{x}\quad\mathcal{L}(x,u,\lambda) = c^{\top}x - \sum_{i=1}^{m}(u_{i}^{\top}A_{i}x+\lambda_{i}b_{i}^{\top}x).
\end{align}

Consider another SOCP problem of the form
\begin{align}
    \minimize{x\in\mathbb{R}^{n}}\quad c^{\top}x,\;
    \st\quad\norm{A_{i}x}^{2}\leqslant b_{i}^{\top}x,\quad\forall i=1,\dots,m,
\end{align}
with the \textit{rotated} second-order cone constraint. To dualize this constraint, we use the fact that for any set of variables $\mu_{i}\in\mathbb{R}^{1}$ , $\lambda_{i}\in\mathbb{R}^{1}$  and $u_{i}\in\mathbb{R}^{n_{i}}$ it holds that 
\begin{align}
&\maximize{u_{i},\mu_{i},\lambda_{i}:\atop\norm{u_{i}}^{2}\leqslant\mu_{i}\lambda_{i}}
-u_{i}^{\top}A_{i}x-\nicefrac{1}{2}\lambda_{i} - \mu_{i}b_{i}^{\top}x 
\nonumber\\
&= 
\maximize{\lambda_{i}\geqslant0}\;-\lambda_{i}(\norm{A_{i}x}^2- b_{i}^{\top}x) 
= 
\left\{
\begin{matrix*}[l]
0, & \text{if}\;\norm{A_{i}x}^2\leqslant b_{i}^{\top}x,\\
-\infty, &\text{otherwise}.
\end{matrix*}
\right.\nonumber
\end{align}
Therefore, the Lagrangian of the SOCP problem writes in variables $x\in\mathbb{R}^{n}$, $\mu,\lambda\in\mathbb{R}^{m}$ and $u\in\mathbb{R}^{n_{i}\times n}$ as 
\begin{align}
  \maximize{\norm{u_{i}}^{2}\leqslant\mu_{i}\lambda_{i}}&\;\minimize{x}\quad\mathcal{L}(x,u,\mu,\lambda) = c^{\top}x \nonumber\\
    &- \sum_{i=1}^{m}(u_{i}^{\top}A_{i}x+\nicefrac{1}{2}\lambda_{i} +\mu_{i}b_{i}^{\top}x).
\end{align}
\end{subequations}

\subsection{Proof of Theorem \ref{th:equilibrium}}\label{proof_th_equilibrium}
Consider the problem of finding an equilibrium solution among the following set of agents. First, consider a price-setter who seeks the optimal prices to coupling constraints \eqref{cc_ref_con_1}--\eqref{cc_ref_con_8} in response to their slacks by solving 
\begin{align}
    &\text{max}_{\lambda^{c},\lambda^{r},\lambda^{w},\lambda^{\varphi},\lambda^{\pi},\lambda^{\underline{\varphi}},\lambda^{\overline{\pi}},\lambda^{\underline{\pi}}}
    \quad\lambda^{c\top}\braceit{
    A\varphi - \vartheta + B\kappa + \delta
    }
    \nonumber\\
    &
    +
    \lambda^{r\top}
    \braceit{
    \mathbb{1} - (\alpha - B\beta)^{\top}\mathbb{1}
    }
    +\lambda^{w\top}\braceit{
    \varphi - \gamma_{1} - \gamma_{2}\pi - \gamma_{3}\kappa
    }
    \nonumber\\
    &
    +\mysum{\ell=1}{E}\lambda_{\ell}^{\varphi}\braceit{
    s_{\ell}^{\varphi} - \norm{F[\grave{\gamma}_{2}(\alpha - \text{diag}[\mathbb{1}]) - \grave{\gamma}_{3}\beta]_{\ell}^{\top}}
    }
    \nonumber\\
    &
    +\mysum{n=1}{N}\lambda_{n}^{\pi}\braceit{
    s_{n}^{\pi} - \norm{F[\breve{\gamma}_{2}(\alpha - \hat{\gamma}_{3}\beta -\text{diag}[\mathbb{1}])]_{n}^{\top}}
    }
    \nonumber\\
    &
    +\mysum{\ell=1}{E}\lambda_{\ell}^{\underline{\varphi}}\braceit{
    \varphi_{\ell} - z_{\hat{\varepsilon}}\norm{F[\grave{\gamma}_{2}(\alpha - \text{diag}[\mathbb{1}]) - \grave{\gamma}_{3}\beta]_{\ell}^{\top}}
    }
    \nonumber\\
    &
    +\mysum{n=1}{N}\lambda_{\pi}^{\overline{\pi}}\braceit{
    \overline{\pi}_{n}-\pi_{n} - z_{\hat{\varepsilon}}\norm{F[\breve{\gamma}_{2}(\alpha - \hat{\gamma}_{3}\beta -\text{diag}[\mathbb{1}])]_{n}^{\top}}
    }
    \nonumber\\
    &
    +\mysum{n=1}{N}\lambda_{\pi}^{\underline{\pi}}\braceit{
    \pi_{n} - \underline{\pi}_{n} - z_{\hat{\varepsilon}}\norm{F[\breve{\gamma}_{2}(\alpha - \hat{\gamma}_{3}\beta -\text{diag}[\mathbb{1}])]_{n}^{\top}}
    }.\label{eq_price_setter_1}
\end{align}
Problem \eqref{eq_price_setter_1} adjusts the prices respecting the slack of each constraint, e.g., $\lambda^{c}\downarrow$ if $ A\varphi - \vartheta + B\kappa > \delta$, and $\lambda^{c}\uparrow$ otherwise. From SOCP property \eqref{eq:dualization_soc}, we know that the last five terms associated with the conic constraints rewrite equivalently as 
\begin{align*}
    &
    -\lambda^{\varphi\top}s^{\varphi}
    -\lambda^{\pi\top}s^{\pi} 
    -\lambda^{\underline{\varphi}\top}\varphi
    -\lambda^{\overline{\pi}\top}\braceit{\overline{\pi}-\pi}
    -\lambda^{\underline{\pi}\top}\braceit{\pi-\underline{\pi}}
    \nonumber\\
    &
    -\mysum{\ell=1}{E}[u^{\varphi} + z_{\hat{\varepsilon}}u^{\underline{\varphi}}]_{\ell}
    F[\grave{\gamma}_{2}(\alpha - \text{diag}[\mathbb{1}]) - \grave{\gamma}_{3}\beta]_{\ell}^{\top}
    \nonumber\\
    &
    - \mysum{n=1}{N}[u^{\pi} + z_{\hat{\varepsilon}}u^{\overline{\pi}} + z_{\hat{\varepsilon}}u^{\underline{\pi}}]_{n}F[\breve{\gamma}_{2}(\alpha - \hat{\gamma}_{3}\beta -\text{diag}[\mathbb{1}])]_{n}^{\top},
\end{align*}
which is linear and separable, and where the dual variables $u^{\varphi},u^{\underline{\varphi}}\in\mathbb{R}^{E\times N}$ and $u^{\pi},u^{\overline{\pi}},u^{\underline{\pi}}\in\mathbb{R}^{N\times N}$ are subject to the following dual feasibility conditions
\begin{subequations}\label{eq:dual_feasibility}
\begin{align}
    &\norm{[u^{\pi}]_{n}}\leqslant \lambda_{n}^{\pi},
    \norm{[u^{\overline{\pi}}]_{n}}\leqslant \lambda_{n}^{\overline{\pi}},
    \norm{[u^{\underline{\pi}}]_{n}}\leqslant \lambda_{n}^{\underline{\pi}},\label{dual_feasibility_1}\\
    &\norm{[u^{\varphi}]_{\ell}}\leqslant \lambda_{\ell}^{\varphi},
    \norm{[u^{\underline{\varphi}}]_{\ell}}\leqslant \lambda_{\ell}^{\varphi},\label{dual_feasibility_2}
    \forall n\in\mathcal{N}, \forall \ell\in\mathcal{E}.
\end{align}
\end{subequations}
By separating the terms with respect to the variables of network assets, network operator, and free terms associated with each consumer, we obtain the revenue functions in \eqref{eq:revenues}. Consider next that each gas supplier $n\in\mathcal{N}$ solves
\begin{subequations}\label{model:profit_gas_injection}
\begin{align}
    \maximize{\vartheta_{n},[\alpha]_{n},c_{n}^{\vartheta},c_{n}^{\alpha}}\quad&\mathcal{R}_{n}^{\text{sup}}(\vartheta_{n},[\alpha]_{n}) - c_{1n}\vartheta_{n} - c_{n}^{\vartheta} - c_{\alpha}^{\vartheta}\label{profit_gas_injection_obj}\\
    \st\quad\lambda_{n}^{\vartheta}\colon
    &\textcolor{white}{z_{\hat{\varepsilon}}}\norm{\grave{c}_{2n}\vartheta_{n}}^{2}\leqslant c_{n}^{\vartheta},\\
    \lambda_{n}^{\alpha}\colon
    &\textcolor{white}{z_{\hat{\varepsilon}}}\norm{F\grave{c}_{2n}[\alpha]_{n}^{\top}}^{2}\leqslant c_{n}^{\alpha},\\
    \lambda_{n}^{\overline{\vartheta}}\colon
    &z_{\hat{\varepsilon}}\norm{F[\alpha]_{n}^{\top}}\leqslant\overline{\vartheta}_{n}-\vartheta_{n},\\
    \lambda_{n}^{\underline{\vartheta}}\colon
    &z_{\hat{\varepsilon}}\norm{F[\alpha]_{n}^{\top}}\leqslant\vartheta_{n}-\underline{\vartheta}_{n},
\end{align}
\end{subequations}
to maximize the profit in response to equilibrium prices. Next, consider that each active pipeline $\ell\in\mathcal{E}$ solves
\begin{subequations}\label{model:profit_active_pipeline}
\begin{align}
    \maximize{\kappa_{\ell},[\beta]_{\ell}}\quad&\mathcal{R}_{\ell}^{\text{act}}(\kappa_{\ell},[\beta]_{\ell})\label{profit_active_pipeline_obj}\\
    \st\quad
    \lambda_{\ell}^{\overline{\kappa}}\colon
    &z_{\hat{\varepsilon}}\norm{F[\beta]_{\ell}^{\top}}\leqslant\overline{\kappa}_{\ell}-\kappa_{\ell},\\
    \lambda_{\ell}^{\underline{\kappa}}\colon
    &z_{\hat{\varepsilon}}\norm{F[\beta]_{\ell}^{\top}}\leqslant\kappa_{\ell}-\underline{\kappa}_{\ell},
\end{align}
\end{subequations}
to maximize the revenue in response to equilibrium prices. Finally, consider a gas network operator which solves
\begin{subequations}\label{model:rent_minimization}
\begin{align}
    \minimize{\pi,\varphi,s^{\pi},s^{\varphi}}\quad&\mathcal{R}^{\text{rent}}(\pi,\varphi,s^{\pi},s^{\varphi})\\
    \st\quad
    \lambda_{r}^{\mathring{\pi}}\colon
    &\pi_{\text{r}} = \mathring{\pi}_{\text{r}}
\end{align}
\end{subequations}
to maximize the network rent in response to equilibrium prices. By taking the path outlined in Appendix \ref{app:dualization}, the first-order optimality conditions (FOC) of equilibrium problems \eqref{model:profit_gas_injection}--\eqref{model:rent_minimization} are given by the following equalities
\begin{subequations}\label{eq:stat_con}
\begin{align}
    \vartheta\colon&
    c_{1} - u^{\vartheta}\circ\grave{c}_{2} - \lambda^{c} + \lambda^{\overline{\vartheta}} - \lambda^{\underline{\vartheta}}  = \mathbb{0},
    \label{eq:stat_con_1}
    \\
    \kappa\colon&
    [\lambda^{c\top}B]^{\top} - [\lambda^{w\top}\gamma_{3}]^{\top} + \lambda^{\overline{\kappa}} - \lambda^{\underline{\kappa}}
    = \mathbb{0},
    \label{eq:stat_con_2}
    \\
    \pi\colon&
    \lambda^{\overline{\pi}} - \lambda^{\underline{\pi}} - [\lambda^{w\top}\gamma_{2}]^{\top} - \mathbb{I}_{r}\circ\lambda^{\mathring{\pi}} = \mathbb{0},
    \label{eq:stat_con_3}
    \\
    \varphi\colon&
    [\lambda^{c\top}A]^{\top} + \lambda^{w} - \lambda^{\underline{\varphi}} = \mathbb{0},
    \label{eq:stat_con_4}
    \\
    s^{\pi}\colon& \lambda^{\pi} = \psi^{\pi},\; s^{\varphi}\colon\lambda^{\varphi} = \psi^{\varphi},
    \label{eq:stat_con_4_a}
    \\
    c^{\vartheta}\colon&\mu^{\vartheta} = \mathbb{1},\; c^{\alpha}\colon\mu^{\alpha} = \mathbb{1},
    \label{eq:stat_con_5}
    \\
    [\alpha]_{n}\colon&
    F
    \braceit{
    u^{\varphi\top}\col{\grave{\gamma}_{2}}_{n} 
    + u^{\pi\top}\col{\breve{\gamma}_{2}}_{n} 
    + z_{\hat{\varepsilon}}[u^{\overline{\vartheta}} + u^{\underline{\vartheta}}]_{n}^{\top}
    }\nonumber\\
    &\quad\quad\quad\quad\quad\quad\!\!\!+F[u^{\alpha}]_{n}^{\top}\grave{c}_{2} 
    + \lambda^{r}= \mathbb{0},
    \label{eq:stat_con_6}
    \\
    [\beta]_{\ell}\colon&
    F
    \braceit{
    u^{\varphi\top}\col{\grave{\gamma}_{3}}_{\ell}
    + u^{\pi\top}\col{\breve{\gamma}_{2}\hat{\gamma}_{3}}_{\ell}
    - z_{\hat{\varepsilon}}[u^{\overline{\kappa}} + u^{\underline{\kappa}}]_{\ell}^{\top}
    }\nonumber\\
    &\quad\quad\quad\quad\quad\;\!+ \mathbb{1}^{\top}\col{B}_{\ell}\lambda^{r} = \mathbb{0}, 
    \label{eq:stat_con_7}
\end{align}
\end{subequations}
where vector $\mathbb{I}_{r}\in\mathbb{R}^{N}$ takes 1 at position corresponding to the reference node, and 0 otherwise. Conditions \eqref{eq:stat_con} are identical to those of centralized problem \eqref{model:cc_ref}, while the set of FOC of problem \eqref{eq_price_setter_1} yields primal constraints \eqref{cc_ref_con_1}--\eqref{cc_ref_con_8}. Together with the primal constraints of problems \eqref{model:profit_gas_injection}--\eqref{model:rent_minimization}, they are identical to the feasibility conditions of the centralized problem. Hence, problem \eqref{model:cc_ref} solves the competitive equilibrium. 

\subsection{Proof of Corollary \ref{cor:revenue_adequcy}}\label{proof:rev_ad}
From the feasibility conditions \eqref{cc_ref_con_1}--\eqref{cc_ref_con_3} and complementarity slackness conditions associted with constraints \eqref{cc_ref_con_4}--\eqref{cc_ref_con_8}, we know that the objective function of the price-setting problem in \eqref{eq_price_setter_1} is zero at optimum. By rearranging the terms of \eqref{eq_price_setter_1}, we have
\begin{align*}
    \mysum{n=1}{N}\mathcal{R}_{n}^{\text{con}} - \mysum{n=1}{N}\mathcal{R}_{n}^{\text{sup}} - \mysum{\ell=1}{E}\mathcal{R}_{\ell}^{\text{act}} = \mathcal{R}^{\text{rent}} + \lambda^{w\top}\gamma_{1}.
\end{align*}
If let $\gamma_{1}=\mathbb{0}$, it remains to show that the congestion rent accumulated by the network is non-negative, i.e., 
\begin{align*}
    &
    \underbrace{
    \braceit{\lambda^{\underline{\varphi}\top} - \lambda^{w\top} - \lambda^{c\top} A}\varphi
    }_{\text{Term A}}
    +
    \underbrace{
    \braceit{\lambda^{w\top}\gamma_{2}+ \lambda^{\underline{\pi}\top}-  \lambda^{\overline{\pi}\top}}\pi
    }_{\text{Term B}}\nonumber\\
    & 
    +
    \underbrace{
    \lambda^{\overline{\pi}\top}\overline{\pi} -\lambda^{\underline{\pi}\top}\underline{\pi}
    }_{\text{Term C}}
    +
    \underbrace{
    \lambda^{\varphi\top} s^{\varphi} 
    +\lambda^{\pi\top}s^{\pi}
    }_{\text{Term D}}
    \geqslant 0.
\end{align*}
From optimality condition \eqref{eq:stat_con_4}, we know that term A is zero. Due to \eqref{eq:stat_con_3}, the term B is zero for all nodes but the reference one, and for the reference node it is $\lambda^{\mathring{\pi}}\mathring{\pi}\geqslant0$ from the dual objective function of problem \eqref{model:rent_minimization}.  Term D is non-negative, because from \eqref{eq:stat_con_4_a} we have that the dual prices $\lambda^{\varphi}$ and $\lambda^{\pi}$ are non-negative, and variables $s^{\varphi}$ and $s^{\pi}$ are lower-bounded by zero as per \eqref{cc_ref_con_4} and \eqref{cc_ref_con_5}. In term C, $\lambda^{\overline{\pi}\top}\overline{\pi}$ and $\lambda^{\underline{\pi}\top}\underline{\pi}$ are non-negative due conditions \eqref{dual_feasibility_1}. Thus, the rent is always non-negative if and only if the network design allows $\underline{\pi}=\mathbb{0}$. 

\subsection{Proof of Corollary \ref{cor:cost_recovery}}
We need to show that the functions \eqref{profit_gas_injection_obj} and \eqref{profit_active_pipeline_obj} are non-negative. Both \eqref{profit_gas_injection_obj} and \eqref{profit_active_pipeline_obj} are lower bounded by their corresponding dual functions, i.e., 
\begin{align*}
    \eqref{profit_gas_injection_obj}\;\geqslant\;&\nicefrac{1}{2}(\lambda_{n}^{\vartheta} + \lambda_{n}^{\alpha}) + \lambda_{n}^{\overline{\vartheta}}\overline{\vartheta}_{n} - \lambda_{n}^{\underline{\vartheta}}\underline{\vartheta}_{n}, &\forall n \in \mathcal{N},\\
    \eqref{profit_active_pipeline_obj}\;\geqslant\;&\lambda_{\ell}^{\overline{\kappa}}\overline{\kappa}_{\ell} - \lambda_{\ell}^{\underline{\kappa}}\underline{\kappa}_{\ell},&\forall \ell \in \mathcal{E}_{a}.
\end{align*}
From the complementary slackness of constraints in \eqref{model:profit_gas_injection} and \eqref{model:profit_active_pipeline}, we know that $\lambda^{\vartheta},\lambda^{\alpha},\lambda^{\overline{\vartheta}},\lambda^{\underline{\vartheta}}\geqslant \mathbb{0}$ and $\lambda^{\overline{\kappa}},\lambda^{\underline{\kappa}}\geqslant \mathbb{0}$. As injection limits are all non-negative, function \eqref{profit_gas_injection_obj} is non-negative if and only if the network design allows $\underline{\vartheta}=\mathbb{0}$. 
As pressure regulation limits for compressors and valves are respectively non-negative and non-positive, function \eqref{profit_active_pipeline_obj} is non-negative if and only if the network design allows $\underline{\kappa}_{\ell}=0,\forall\ell\in\mathcal{E}_{c}$ and $\overline{\kappa}_{\ell}=0,\forall\ell\in\mathcal{E}_{v}$.

\nomenclature[A]{$\mathcal{N}$}{Set of nodes.}
\nomenclature[A]{$\mathcal{E}$}{Set of pipelines.}
\nomenclature[A]{$\mathcal{E}_{a},\mathcal{E}_{c},\mathcal{E}_{v}$}{Set of active, compressor-, valve-hosting pipelines.}

\nomenclature[B]{$p$}{Probability of violating performance guarantee.}
\nomenclature[B]{$v$}{Confidence level of performance guarantee.}
\nomenclature[B]{$A$}{Node-edge incidence matrix.}
\nomenclature[B]{$\delta$}{Vector of nominal gas extraction rates.}
\nomenclature[B]{$w$}{Vector of pipeline friction coefficients.}
\nomenclature[B]{$\underline{\rho},\overline{\rho}$}{Vectors of min. and max. pressure limits.}
\nomenclature[B]{$\underline{\pi},\overline{\pi}$}{Vectors of min. and max. squared pressure limits.}
\nomenclature[B]{$\underline{\kappa},\overline{\kappa}$}{Vectors of min. and max. squared regulation limits.}
\nomenclature[B]{$B$}{Sending node - active pipeline incidence matrix.}
\nomenclature[B]{$b$}{Vector of gas mass - pressure conversion factors.}
\nomenclature[B]{$c_{1},c_{2}$}{Vectors of the 1\textsuperscript{st}- and 2\textsuperscript{nd}-order cost coefficients.}
\nomenclature[B]{$\grave{c}_{2}$}{Factorization of the 2\textsuperscript{nd}-order cost coefficients.}
\nomenclature[B]{$\underline{\vartheta},\overline{\vartheta}$}{Vectors of min. and max. gas injection limits.}
\nomenclature[B]{$\gamma,\hat{\gamma},\breve{\gamma},\grave{\gamma}$}{Linearization coefficients (and their transformations) associated with the Weymouth equation.}
\nomenclature[B]{$\Sigma,F$}{Forecast errors covariance matrix and its factorization.}
\nomenclature[B]{$\varepsilon,\hat{\varepsilon}$}{Joint and individual constraint violation parameters.}
\nomenclature[B]{$z_{\hat{\varepsilon}}$}{Safety parameter in chance constraint reformulation.}
\nomenclature[B]{$\psi^{\varphi}$}{Vector of flow variance penalty factors.}
\nomenclature[B]{$\psi^{\pi}$}{Vector of pressure variance penalty factors.}
\nomenclature[B]{$\mathcal{R}^{(\cdot)}$}{Revenue associated with network agent $(\cdot)$.}
\nomenclature[B]{$S$}{Sample complexity in out-of-sample analysis.}

\nomenclature[C]{$\vartheta$}{Vector of nodal gas injections.}
\nomenclature[C]{$\kappa$}{Vector of pressure regulation rates.}
\nomenclature[C]{$\pi$}{Vector of nodal pressures.}
\nomenclature[C]{$\varphi$}{Vector of gas flows.}
\nomenclature[C]{$s^{\pi}$}{Vector of pressure standard deviations.}
\nomenclature[C]{$s^{\varphi}$}{Vector of flow standard deviations.}
\nomenclature[C]{$c^{\vartheta},c^{\alpha}$}{Vectors of 2\textsuperscript{nd}-order nominal and recourse costs.}
\nomenclature[C]{$\alpha$}{Matrix of gas injection recourse decisions.}
\nomenclature[C]{$\beta$}{Matrix of pressure regulation recourse decisions.}
\nomenclature[C]{$\lambda,u$}{Vectors and matrices of dual variables.}

\printnomenclature

%\balance
\bibliographystyle{ieeetr}
\bibliography{references}

\begin{IEEEbiography}[{\includegraphics[width=1in,height=1.25in,clip,keepaspectratio]{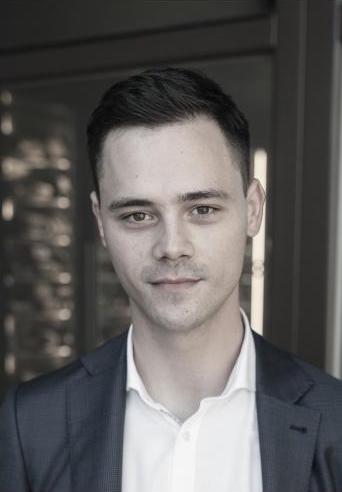}}]{Vladimir Dvorkin} (S’18, M'21) received a Ph.D. degree from the Technical University of Denmark, Lyngby, Denmark, in 2021. He is currently a postdoctoral associate with MIT Energy Initiative and Laboratory for Information and Decision Systems of Massachusetts Institute of Technology, Cambridge, MA, USA. His research interests include economics, game theory, optimization, and their applications to energy systems and markets. 
\end{IEEEbiography}
\begin{IEEEbiography}[{\includegraphics[width=1in,height=1.25in,clip,keepaspectratio]{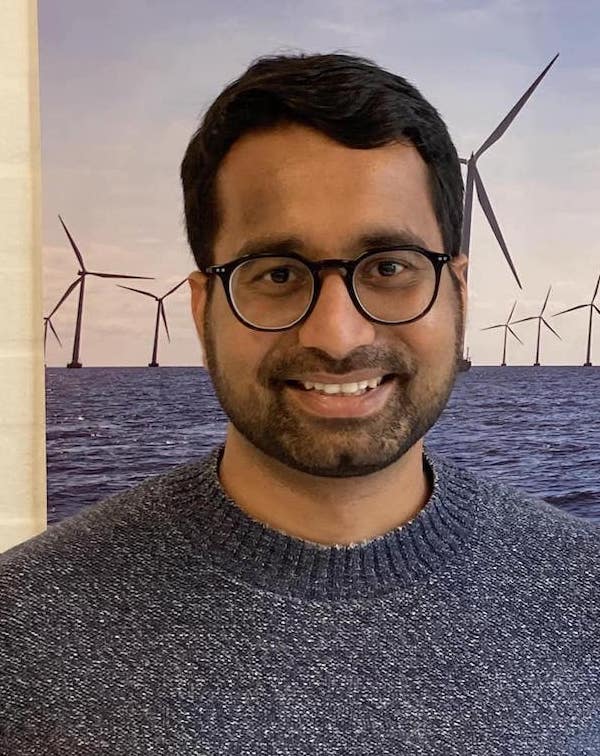}}]{Anubhav Ratha} (S’18), received his M.Sc. degree in energy science and technology from ETH Zurich and is currently pursuing his PhD with the Department of Electrical Engineering, Technical University of Denmark. His research focus is on the design of market mechanisms and products for integrated energy systems of the future. Before his PhD studies, he co-founded a startup leveraging behavioral demand response to develop solutions for the green transition of electricity systems.
\end{IEEEbiography}
\begin{IEEEbiography}[{\includegraphics[width=1in,height=1.25in,clip,keepaspectratio]{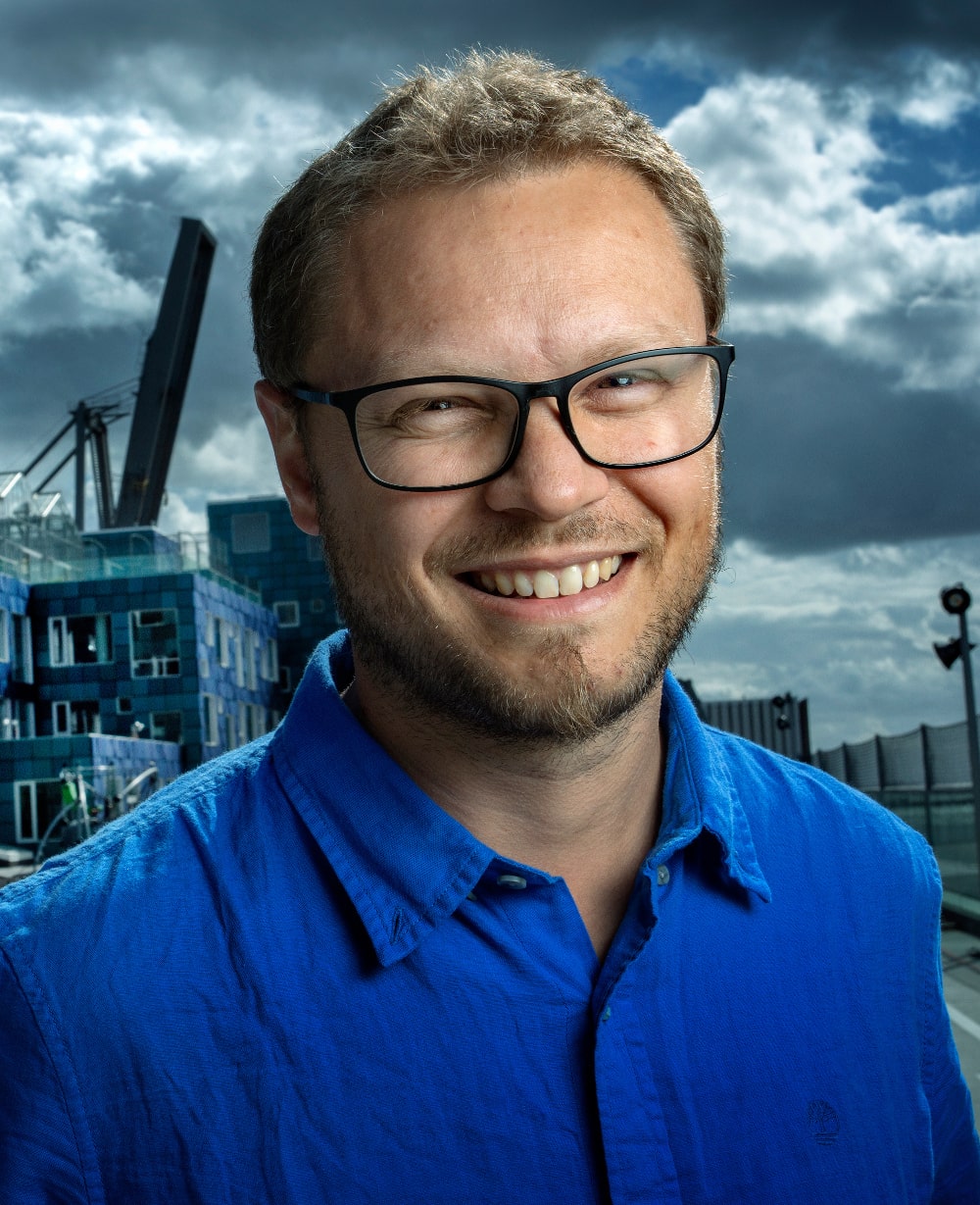}}]{Pierre Pinson} (SM'13, F'20) received his M.Sc. degree in applied mathematics from the National Institute for Applied Sciences, Toulouse, France, and the Ph.D. degree in energetics from Ecole des Mines de Paris, Paris, France. He is a Professor of Operations Research at the Department of Technology, Management and Economics of the Technical University of Denmark, Lyngby, Denmark. His research interests include forecasting, uncertainty estimation, optimization under uncertainty, decision sciences, and renewable energies. He is the Editor-in-Chief for the International Journal of Forecasting.
\end{IEEEbiography}
\begin{IEEEbiography}[{\includegraphics[width=1in,height=1.25in,clip,keepaspectratio]{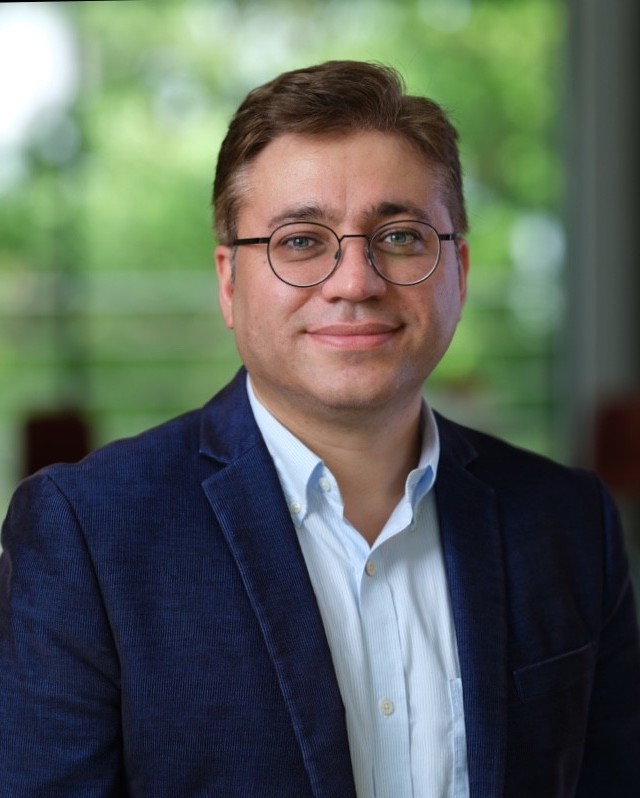}}]{Jalal Kazempour} (SM'18) received  the Ph.D. degree in electrical engineering from the University of Castilla-La Mancha, Ciudad Real, Spain, in 2013. He is currently an Associate Professor with the Department of Electrical Engineering, Technical University of Denmark, Lyngby, Denmark. His focus area is the intersection of multiple fields, including power and energy systems, electricity markets, optimization, and game theory.
\end{IEEEbiography}

\endgroup
\end{document}